\documentclass{article}

\usepackage{amssymb,graphicx,epsf}
\makeatother

\begin{document}

\date{}

\title{\Large\bf Nodal Parity Invariants of Knotted Rigid Vertex Graphs}
\author{Louis
H. Kauffman\\ Department of Mathematics, Statistics \\ and Computer Science (m/c
249)    \\ 851 South Morgan Street   \\ University of Illinois at Chicago\\
Chicago, Illinois 60607-7045\\ $<$kauffman@uic.edu$>$\\ and \\ Rama Mishra
\\ IISER, Pune, India\\
$<$r.mishra@iiserpune.ac.in$>$}

\maketitle

\subsection*{\centering {Abstract}}

{\em This paper introduces new invariants of rigid vertex graph embeddings by using 
non-local combinatorial information that is available at each graphical node.}

\section{Introduction}

\noindent A rigid vertex 4-valent graph (briefly, an RV4 graph)is a 4-valent graph whose vertices are replaced by rigid 2-disks or 3-balls. Each disk or ball has four strands attached to it. By a  (possibly) knotted RV4 graph we mean an embedding of a RV4 graph into $\mathbb{R}^3.$ Such a graph is said to be
{\em unknotted} if the embedding is rigid vertex isotopic (see the definition below) to an embedding of 
the graph in the plane. In a knotted RV4 graph a rigid vertex typically appears as in Figure 1. 

     \begin{center}
     \begin{tabular}{c}
     \includegraphics[height=2cm]{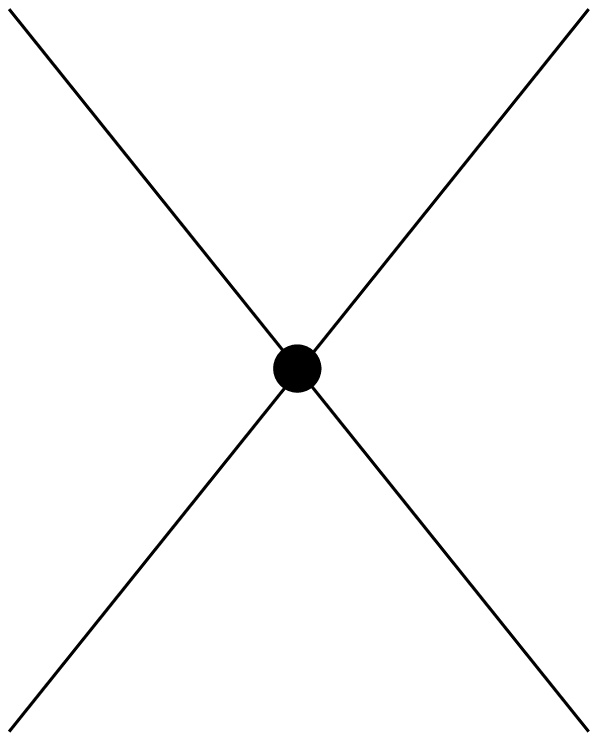}
     \end{tabular}
     \end{center}

\begin{center}

{\bf Figure 1 - A Rigid Vertex}

\end{center}

\noindent In the case of an oriented RV4 graphs, there are  two possible orientations at a given node,  as shown in Figure 2.

     \begin{center}
     \begin{tabular}{c}
     \includegraphics[height=4cm]{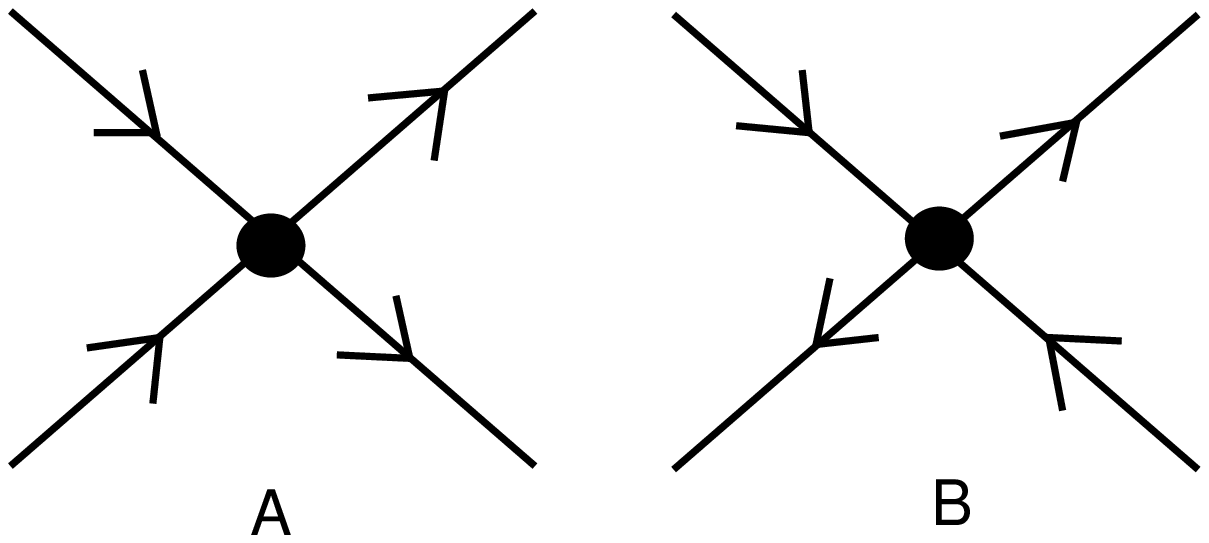}
     \end{tabular}
     \end{center}

\centerline{\bf Figure 2 - Oriented Vertices.}
\vspace{.1in}
   
\noindent Two embedded (possibly knotted) RV4 graphs are RV equivalent (or RV ambient isotopic or just isotopic) if their graph diagrams are transformed to each other by a finite sequence of generalized Reidemeister moves \cite{kau1} (I), (II), (III) as shown in Figure 3 and (IV) and (V)  as shown 
in Figure 3.1.

\begin{center}
 \begin{tabular}{c}
\includegraphics[height=5cm]{ReidMoves.EPSF}
 \end{tabular}
\end{center}

\centerline{\bf Figure 3 - Reidemeister Moves.}

\vspace{.1in}
     
     \begin{center}
 \begin{tabular}{c}
\includegraphics[height=8cm]{RVMoves.EPSF}
 \end{tabular}
\end{center}

\centerline{\bf Figure 3.1 - Rigid Graph Moves (IV) and (V).}

\vspace{.1in}

\noindent It is clear that in a rigid vertex graph $G$, if we replace each rigid vertex (node) by any tangle $T$ we will obtain a link diagram $T_G$. It is easy to see that if two rigid vertex graphs $G_1$ and $G_2$ are isotopic, then the corresponding link diagrams $T_{G_1}$ and $T_{G_2}$ will be isotopic for every choice of $T$. For example we can replace each node of an RV graph by the tangle $T$ as shown in Fig 4. Here the tangle $T$ is a single positive crossing.
\bigbreak

We shall obtain link diagrams whose isotopy classes give  invariants for rigid vertex graphs. Here is  an example. In Figure 5, we have two rigid vertex graphs $G$ and $H$. If we make $T_G$ and $T_H$ we see in Figure 5 that the resulting graphs are not isotopic. Hence the graphs $G$ and $H$ are not isotopic. Note that the local replacements for $T_G$ and $T_H$ are identical positive crossings.

\begin{center}
 \begin{tabular}{c}
\includegraphics[height=2cm]{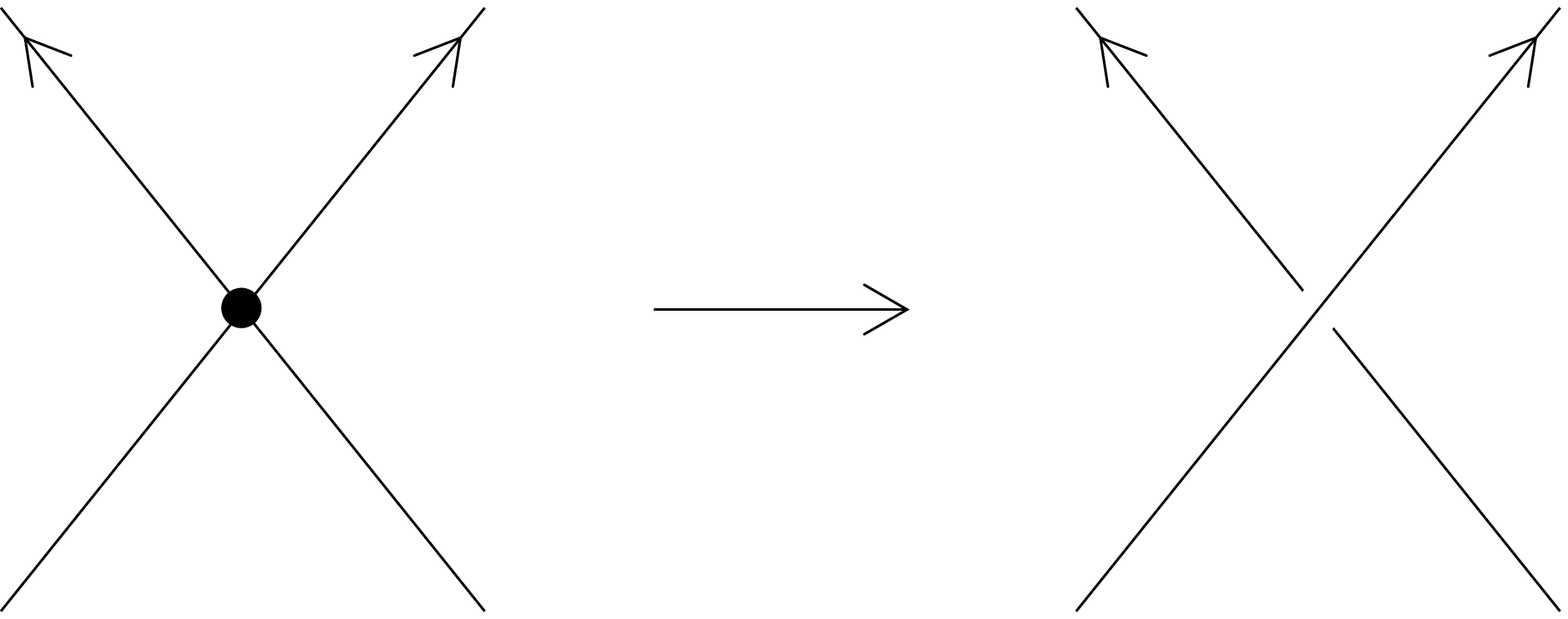}
 \end{tabular}
\end{center}

\centerline{\bf Figure 4 - Positive Replacement.}

\begin{center}
 \begin{tabular}{c}
\includegraphics[height=4cm]{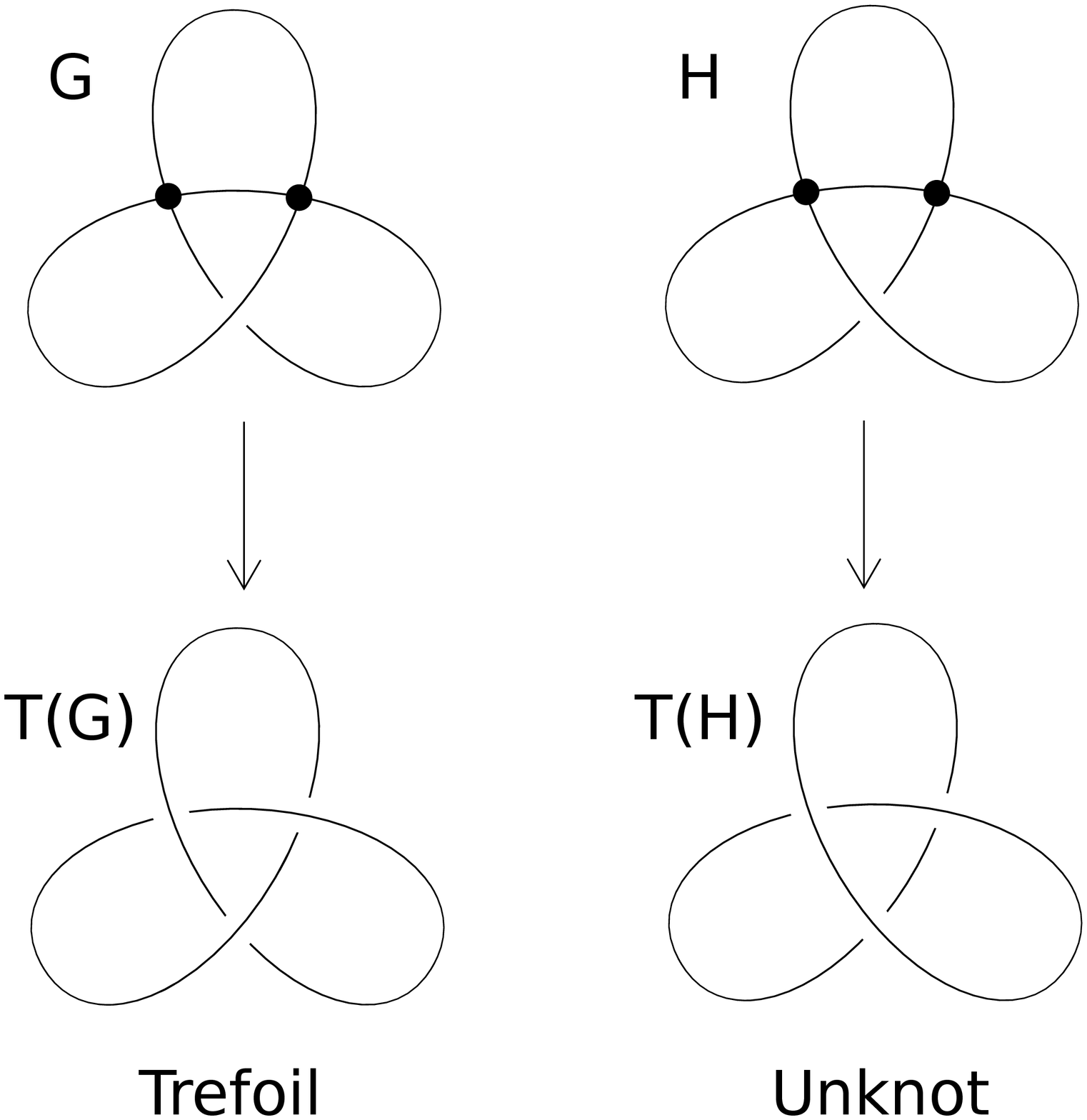}
 \end{tabular}
\end{center}

\centerline{\bf Figure 5 - Positive Replacements.}

\vspace{.1in}

If the rigid vertex graphs have more nodes, and if we replace each node by the same tangle $T$ we may not be able to distinguish the graphs.
For example in Figure 6, if we replace all three nodes by a negative crossing, we obtain an unknotted diagram. However if we replace each node by a positive crossing then the resulting diagram is a right handed trefoil.

\begin{center}
 \begin{tabular}{c}
\includegraphics[height=6cm]{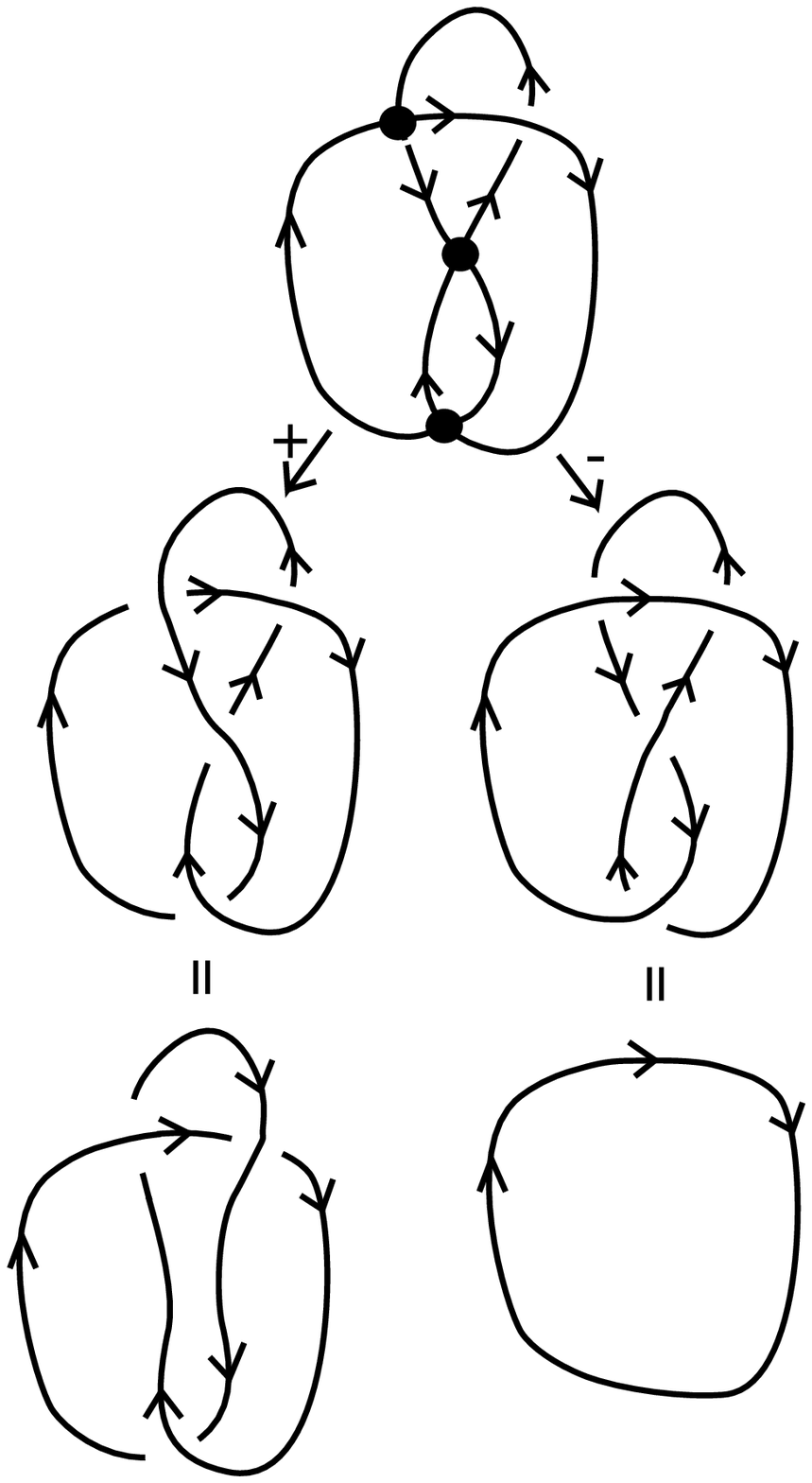}
 \end{tabular}
\end{center}

\centerline{\bf Figure 6 - Positive and Negative Replacements.}

\vspace{.1in}

Thus, replacing each node by the same tangle does not give us full information about the rigid vertex graph. In \cite{kau1} it was shown that the ambient isotopy class of a collection of links $C(\cal{G})$ associated to an RV4 graph $\cal{G}$ (with nodes only of type A or else just an unoriented graph) obtained by making $4$ possible replacements at each node for unoriented graphs (two crossing choices and two smoothing choices) and $3$ possible orientation choices for graphs with type
A orientations (two crossing choices and one oriented smoothing choice)  is an invariant of the RV equivalence of the graph $\cal{G}$. This collection of ambient isotopy types $C(G)$ is quite a strong invariant, however for graphs having large number of nodes it is very hard to determine the set $C(G)$. Thus, it is desirable to look for invariants for RV graphs which are easy to find and can detect the properties of the graphs. If one could categorize the vertices using some property which remains invariant under all Reidemeister Moves, then we can use different replacements for each type of node. Thus we will obtain one link diagram associated to a rigid vertex graph and the isotopy class of link thus associated will be an invariant for RV graphs. The simplest way to execute this idea is to define the notion of parity (even/odd) for each node and replace each even node by a positive crossing and each odd node by a negative crossing (or vice-versa). We will see that the link obtained is an invariant. Thus, the main thing is to have a definition of parity for rigid vertex graphs. We will use Gauss codes to define parity for connected RV graphs. Here we use connectivity in the usual graph theoretic sense: {\em A graph is connected if, given two vertices of the graph there is a path in the graph from one to the other.}
\smallbreak

On the other hand, we wish to speak of the {\it link components} of a rigid vertex graph $G.$ These components are not topological components of the graph. They correspond to traverses of the graph that
utilize the cyclic order that is given at each rigid vertex. A rigid vertex or a diagrammatic crossing has a cyclic order of its four edges dicatated by its embedding in the plane that is given, lets say by $[a,b,c,d].$ We then say that the edges $a$ and $c$ are {\it opposite} and that the edges $b$ and $d$ are opposite. A {\it link component} of $G$ consists in a subgraph obtained by taking a walk on $G$ so that at every crossing and every rigid vertex node one goes in by one edge and out by an opposite edge. Note that the walk may use all four edges of a crossing or a node, or it may use only two of them.
\smallbreak

{\it For graphs with more than one link component, we differentiate vertices belonging to just one component and the vertices that are common between two components and define them to be even or odd accordingly.} We call this the {\it link parity convention.} Note that the link parity convention is very simple and just distinguishes between self-crossings and inter-component crossings. In the next section we will define a subtler parity for self-crossing nodes of a graph.  We will see that a link-type  thus obtained by resolving the nodes according to parity is an invariant. We will see in examples that making different replacements at the nodes of different parity provides better invariants.
\bigbreak

This paper is organized as follows: In section 2 we define the notion of parity at every node of a rigid vertex graph . In section 3 the definition of the invariant is given that depends upon parity and/or distance. Section 4 deals with examples. In section 5 we introduce the concept of special oriented nodes and special oriented nodes with indicators. Section 6 deals with formulating RNA folding as a rigid vertex graphs with or without indicators, and show how these parity invariants are useful in classifying foldings. In section 7 we discuss virtual knot theory in brief and explore how these parity invariants can be used when in a rigid vertex graph there are virtual nodes as well as graphical nodes. 
\bigbreak

The methods for making invariants of rigid vertex graphs that we use in this paper originate from the papers \cite{kau1,kau2}. The concepts of using parity for graphs that we are utilizing here originate
in the paper \cite{SL} by the first author and in the work of Vassily Manturov \cite{MP1,MP2}. In both
of these previous cases, parity is used in working with virtual knot theory. 
\bigbreak

\section {The Concept of Parity for Rigid Vertex graphs}

If the nodes in the rigid vertex graphs can be grouped into two distinct types on the basis of some property which is preserved under Reidemister Moves then one can define new invariants using this information. This property is known as {\em parity} and the two types are called {\em even} and {\em odd}. Parity can be defined in different ways in different situations.  

\subsection{Parity in a single component RV graph}

We define the notion of Gauss code for oriented knotted RV graphs, all of whose nodes are of type $A$ as in Figure $2$, in the same manner as it is done for the oriented knot diagrams. We pick a point on the RV graph that is not a node and then move along the graph following the given orientation and proceeding along the opposite edge to a given incoming edge at the node. We start labeling the distinct nodes by distinct numerals $1,2,3,\ldots$ and when we meet the same node second time we label it again by the same numeral that we used the first time. Thus each numeral is used precisely twice. For example the Gauss code for the RV graph shown in Figure 7 will be $123132$. Thus each node in the graph is identified with a label which is a natural number. For a particular node labeled by $k$ (i.e, the number $k$ occurs twice in the Gauss code) the number of labels that occur in between the two occurrences of $k$ is called the {\it nodal distance} of the node $k$. If the nodal distance of a vertex is odd we say that the vertex has { \it odd parity } otherwise it is of {\it even parity.} In the above example, the nodes $2$ and $3$ are of odd parity while the node $1$ has even parity. It is easy to see that the parity of each node in a RV graph remains invariant under all the Reidemeister Moves for RV graphs.

\begin{center}
\begin{tabular}{c}
\includegraphics[height=4cm]{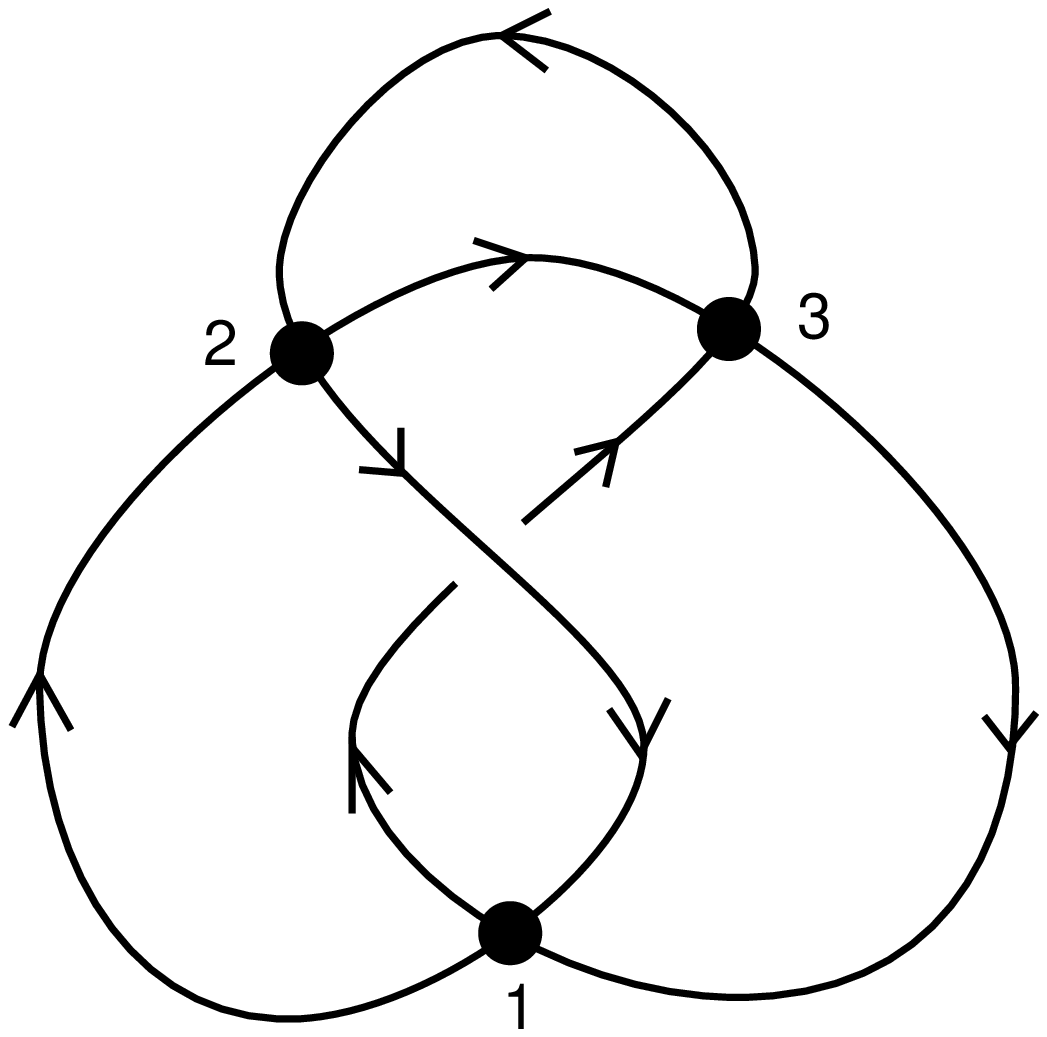}
\end{tabular}
     
\end{center}
     
\centerline{\bf Figure 7 - Gauss Code $123132.$}

\subsection{Parity in a RV graph with multiple link components}  
 
In a RV graph with more than one link component, the nodes may be of two types. In one type the node lies on only one of the link components and the other type it is common between two link components. We call the former as an odd node and the later as an even node. For example in Figure 8, the node 2 is even while the nodes 1 and 3 are odd. This is the {\it link parity convention.}.

\vspace{-.1in}

\begin{center}
     \begin{tabular}{c}
     \includegraphics[height=2cm]{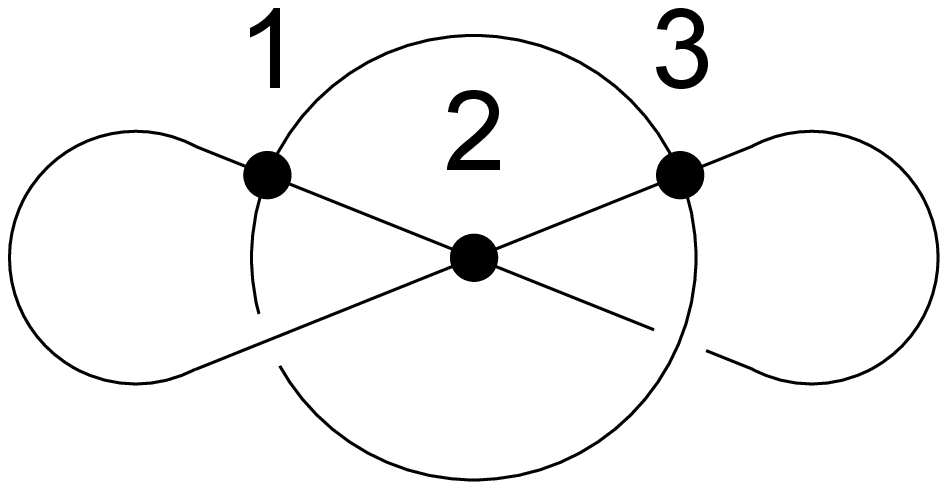}
     \end{tabular}
     \end{center}
     
    \centerline{\bf Figure 8 - Link Parity.}

\section{ Parity invariant for knotted  RV graphs}

It was discussed in \cite{kau1} and \cite{kau2} that if we replace each node of a RV graph $G$ with a tangle to form a link $L(G)$ then any rigid equivalence of $G$ induces a corresponding equivalence of $L(G)$. We have noted that the nodal parity of a RV graph $G$ remains unchanged in any RV isotopy. Thus we can replace nodes of different parity by different tangles. Let us make the simplest choice: replace each node with even parity with a positive crossing  and each node with odd parity with a negative crossing.  This will transform the RV graph $G$ into a link $L(G)$. The ambient isotopy class of $L(G)$ is an RV isotopy invariant for $G.$  We see that $L(G)$ provides us more information than $T(G)$. In Figure 9, we see that the Gauss code for the graph is $123132$ and thus the node $1$ is even, whereas the vertices $2$ and $3$ are odd. So if we replace the even node by a positive crossing and the odd nodes by a negative crossing we obtain the figure eight knot (Figure 9).

\begin{center}
     \begin{tabular}{c}
     \includegraphics[height=4cm]{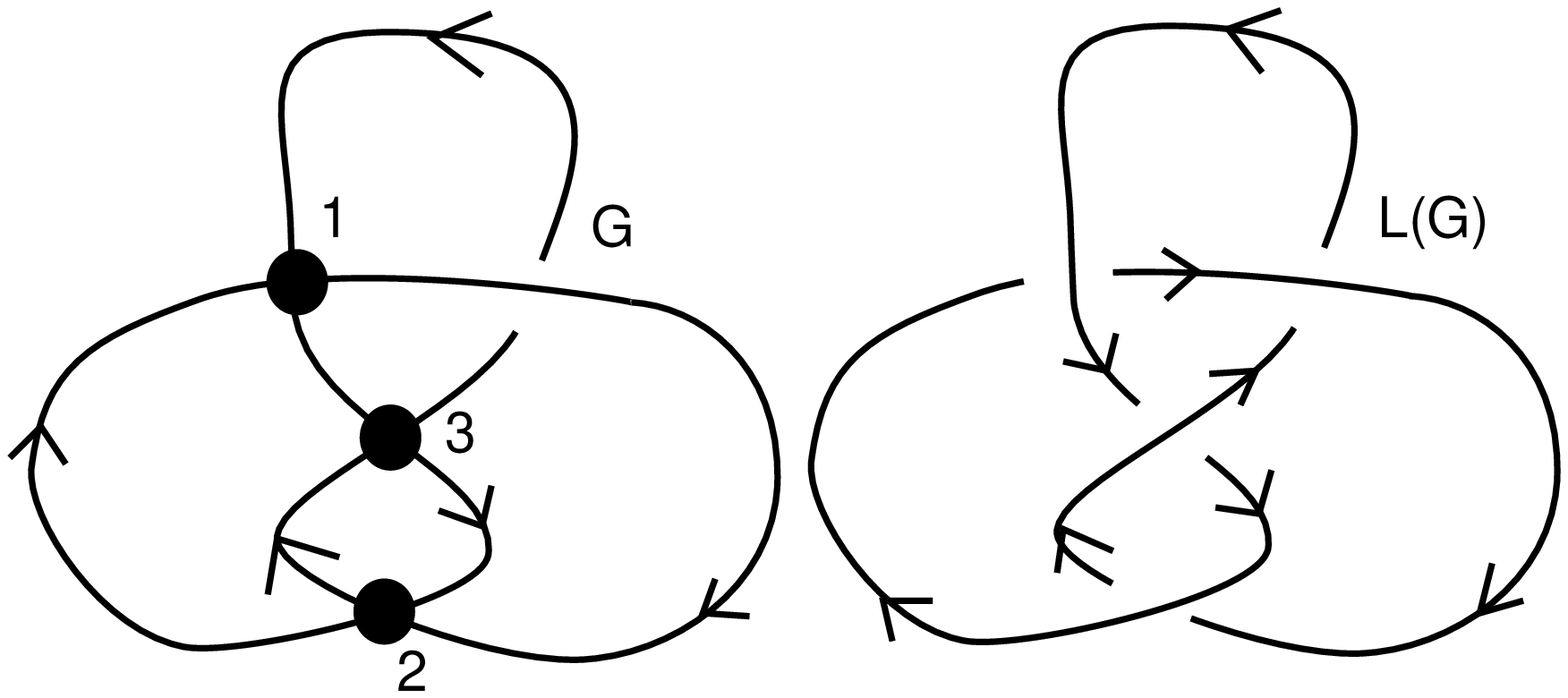}
     \end{tabular}
     \end{center}
     
    \centerline{\bf Figure 9 - Parity Replacement.}

\section{ Examples}

\begin{enumerate}

\item  Let us consider two oriented RV graphs $G$ and $H$ shown  Figure 10.  We note that $G$ and $H$ are mirror images of each other. The Gauss code for both the graphs is $1212$. Thus both the nodes are of odd parity. Replacing both nodes by a negative crossing the graph $G$ results in the trefoil knot whereas the graph $H$ transforms into an unknot. This shows that the graphs $G$ and $H$ are not rigid vertex isotopic, giving the same result as we obtained in Figure 5.

 \item  The RV graphs $G_1$ and $G_2$ shown in Figure 11 are not RV isotopic as the knots $L(G_1)$ and $L(G_2)$ obtained are not ambient isotopic.

    \begin{center}
     \begin{tabular}{c}
     \includegraphics[height=4cm]{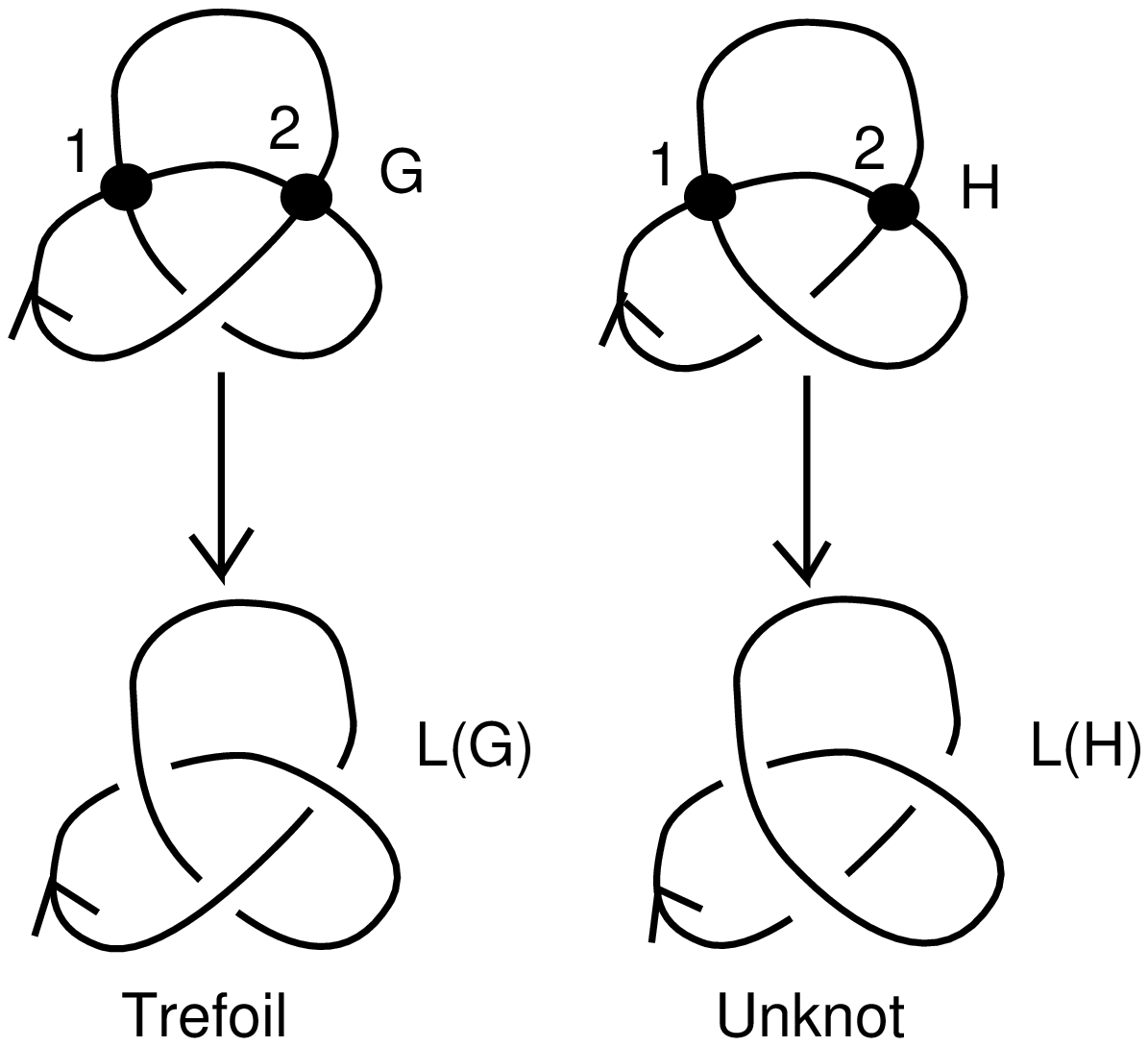}
     \end{tabular}
     \end{center}
     
    \centerline{\bf Figure 10 - Parity Replacement.}

 \begin{center}
     \begin{tabular}{c}
     \includegraphics[height=7cm]{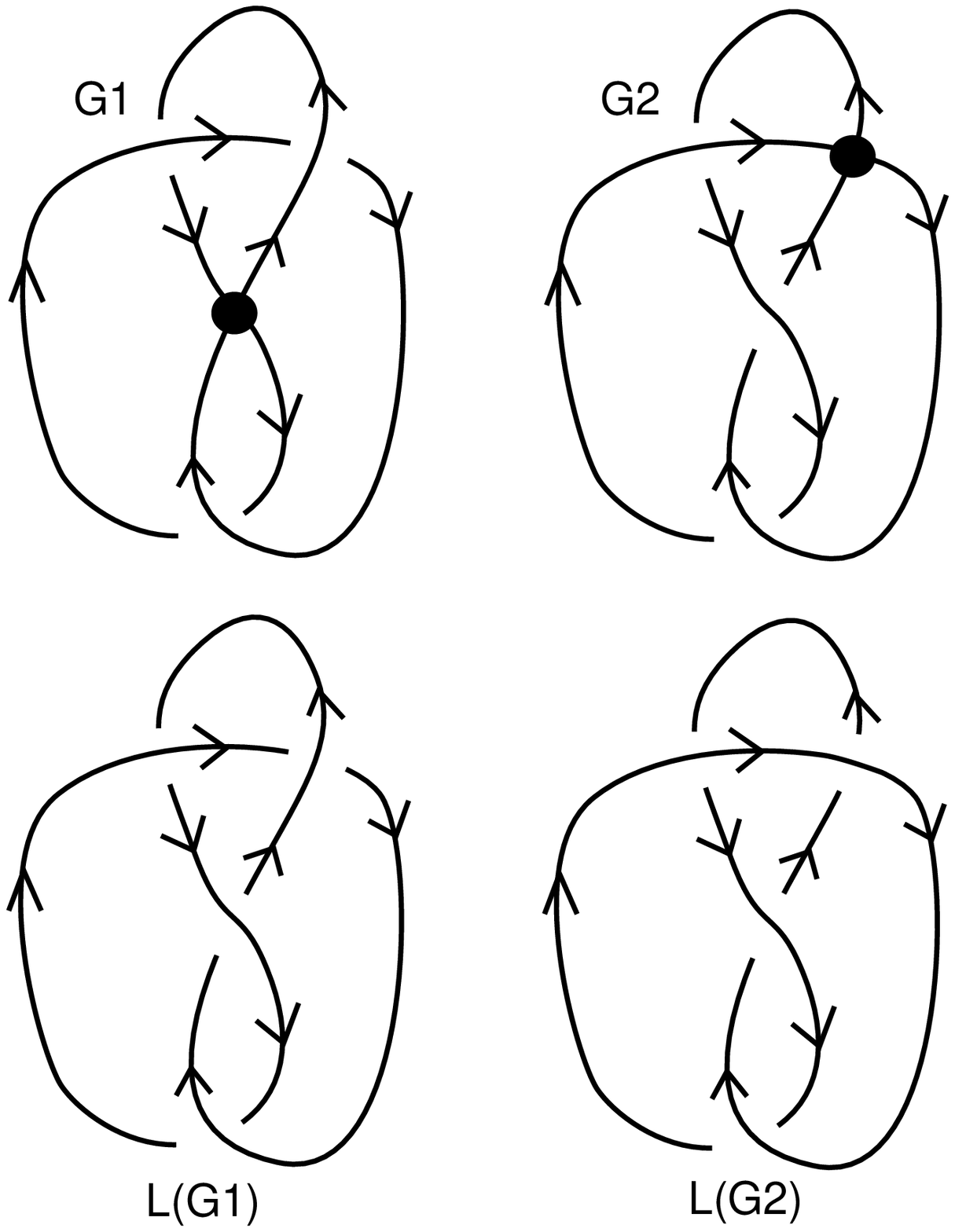}
     \end{tabular}
     \end{center}
     
    \centerline{\bf Figure 11 - Parity Replacement.}

\item   Let $G$ be a rigid vertex graph as shown in Figure 12. Let $G^{*}$ denote its mirror image.

\begin{center}
     \begin{tabular}{c}
     \includegraphics[height=10cm]{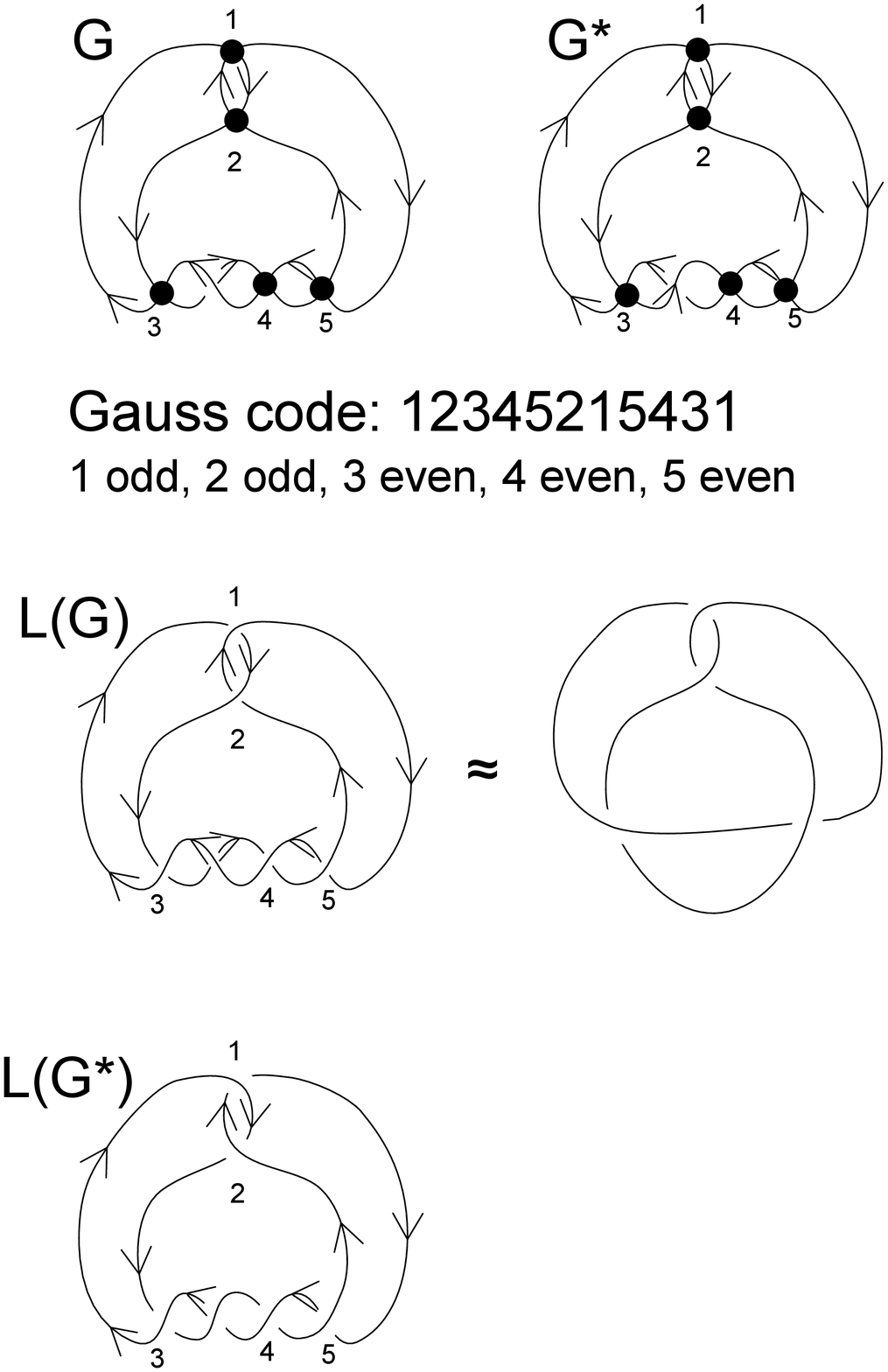}
     \end{tabular}
     \end{center}
     
    \centerline{\bf Figure 12 - Distinct Mirror Images.}

\item  Let $G$ be a graph as shown in Figure 13.   The mirror image of this graph is
$G^{*}$ as shown in the figure. Upon taking parity replacement, we obtain the knots
$L(G)$ and $L(G^{*}.$ A calculation of their Jones polynomials (that we omit here) reveals that the two 
knots are distinct. This proves that the original rigid vertex graph $G$ is chiral. Note that here we have
used parity to determine chirality. 
    
    \begin{center}
     \begin{tabular}{c}
     \includegraphics[height=6cm]{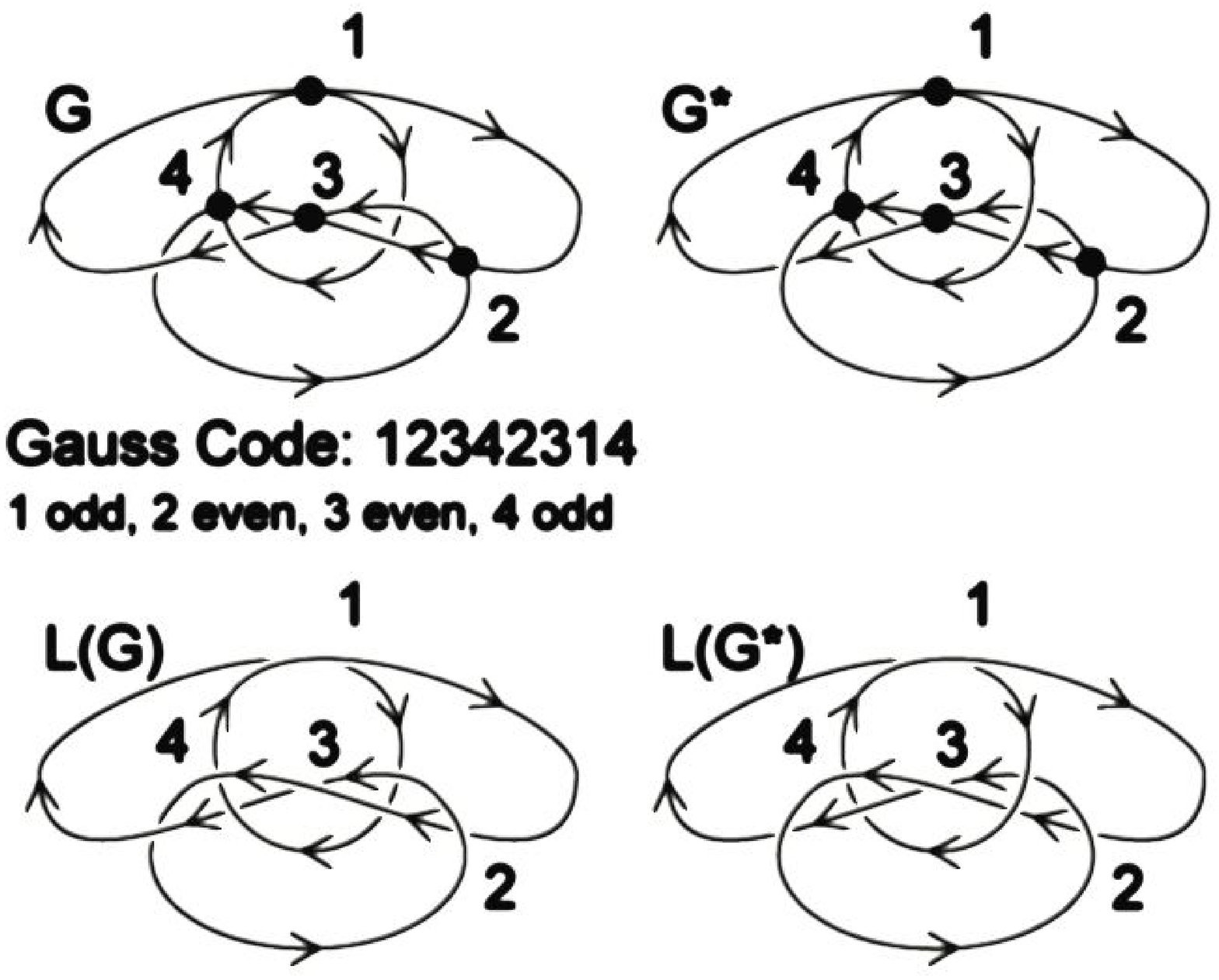}
     \end{tabular}
     \end{center}
     
    \centerline{\bf Figure 13 - Distinct Mirror Images.}

\item Consider the rigid vertex graphs $G1$ and $G2$ as shown in Figure 14. Here $L(G2)$ is isotopic to a trefoil knot while  $L(G1)$ is isotopic to the trivial knot and hence $G1$ and $G2$ are not RV isotopic. Note that the underlying graphs have the same Gauss code.

\item Consider $G_1$ and $G_2$ as shown in Figure 15. These are both RV graphs with two link components. In both the graphs we have one node with even parity and one with odd parity. Here we are using the link parity convention of section 2.2.
Thus if we replace an even node with a positive crossing and the odd node with a negative crossing we see that the link diagram thus obtained are not isotopic. Hence the graphs $G1$ and $G2$ are not isotopic.

\item  Consider an example of a rigid vertex graph $G$ with 3 link components as shown in Figure 16. We can clearly see that $G$ and its mirror image $G^{*}$ do not resolve to the same link. In fact, here
we have $L(G)$ is unlinked while $L(G^{*})$ is the non-trivial linking of the Borommean rings.

 \begin{center}
     \begin{tabular}{c}
     \includegraphics[height=7cm]{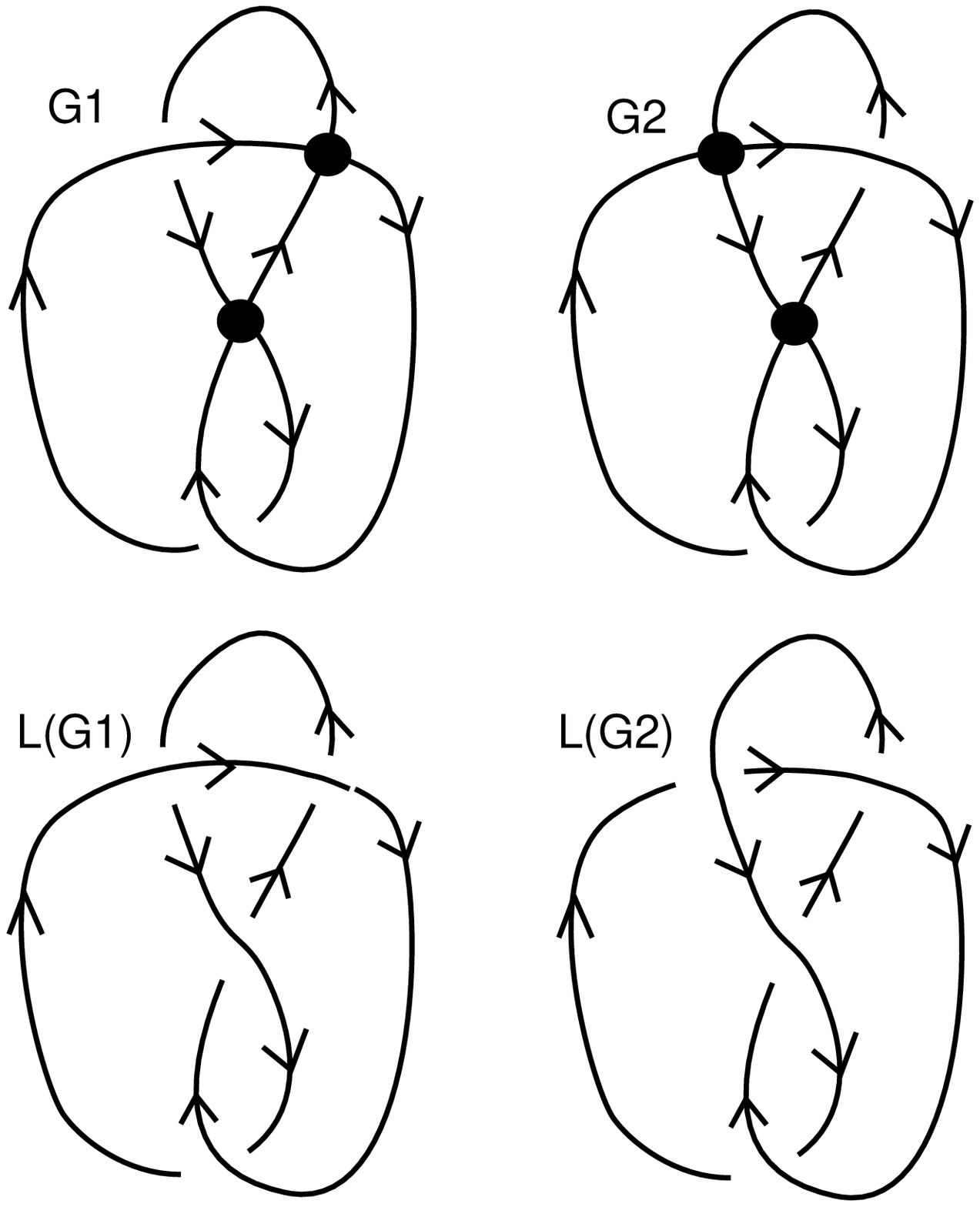}
     \end{tabular}
     \end{center}
     
    \centerline{\bf Figure 14 - Shifting a Node.}
    
\end{enumerate}

\noindent {\bf\large Remark} We have observed in these examples above that in a rigid vertex graph resolving each node by some tangle gives us an invariant of rigid vertex isotopy. We also observed that if we could use different tangles to replace two vertices of different parity we obtain more refined invariant. If we can categorize vertices in a number of groups and use different replacements for each one of them, we may obtain more refined invariants.This sort of situation arises in the folding of RNA molecules, see section 6. 

\begin{center}
     \begin{tabular}{c}
     \includegraphics[height=5cm]{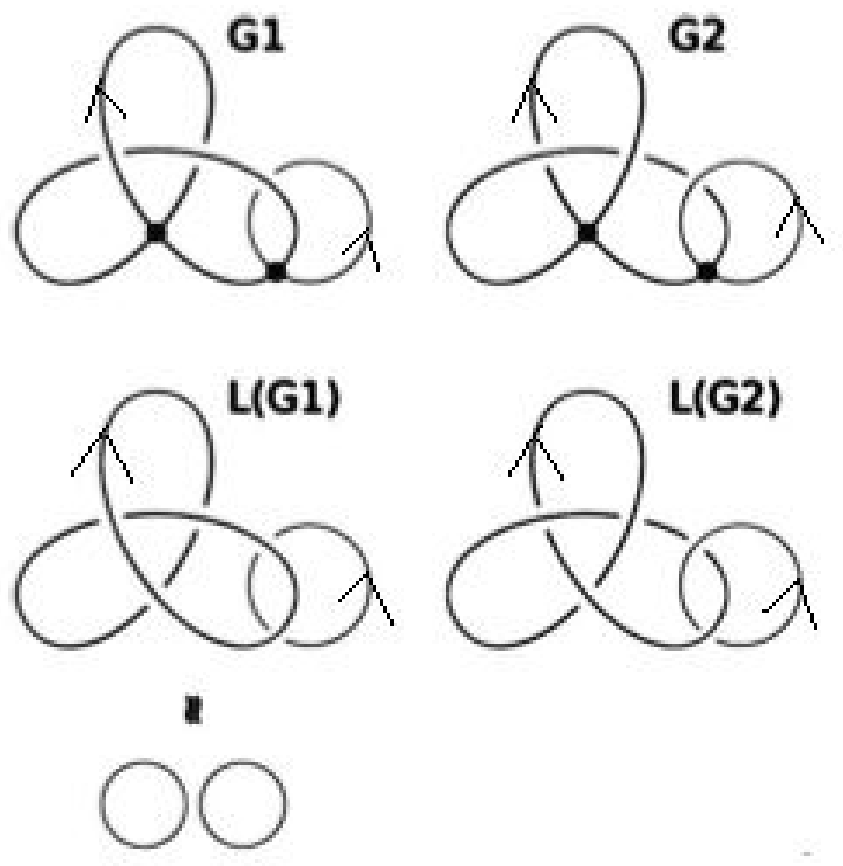}
     \end{tabular}
     \end{center}
     
    \centerline{\bf Figure 15 - Parity Replacement in Multiple Link Components.}

\begin{center}
     \begin{tabular}{c}
     \includegraphics[height=6cm]{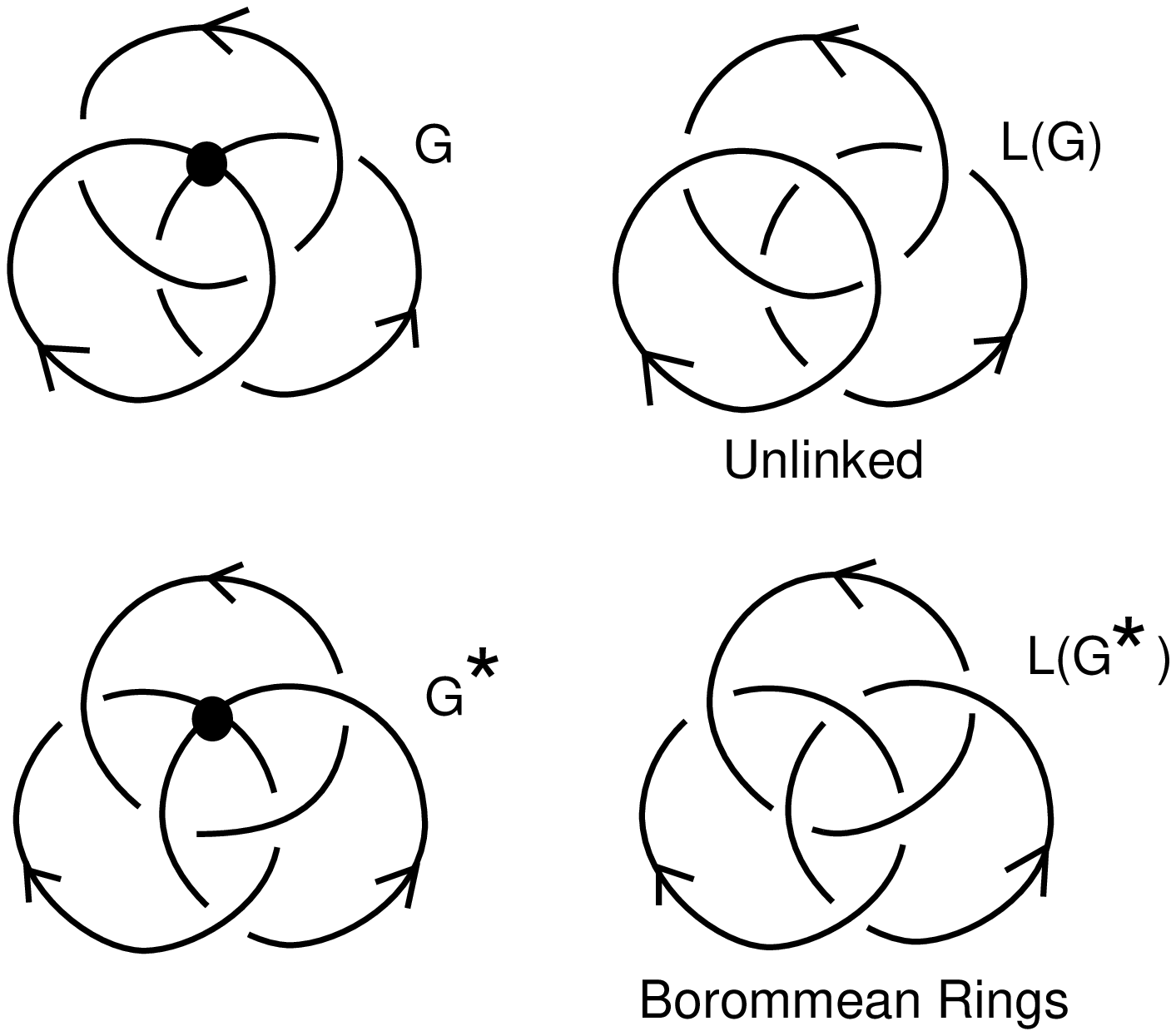}
     \end{tabular}
     \end{center}
     
    \centerline{\bf Figure 16 - Mirror Distinction.}

\section{Rigid vertex graph with special oriented node}

Suppose we have oriented rigid vertex graphs. Then we have two types of nodes as shown in Figure 2 (A and B). The node in Figure 2 (A) is a usual oriented node whereas the node in Figure 2 (B) is a special oriented node. We can think of a special oriented node as in Figure 17.

\begin{center}
     \begin{tabular}{c}
     \includegraphics[height=2cm]{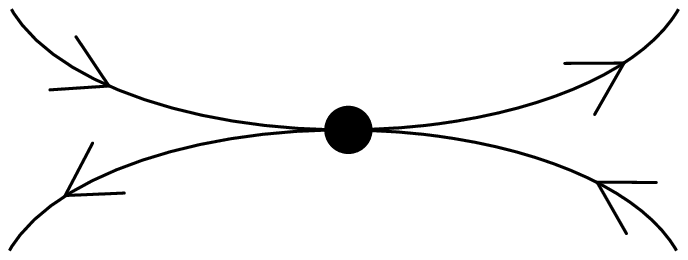}
     \end{tabular}
     \end{center}
     
    \centerline{\bf Figure 17 - Special Oriented Node.}

\vspace{.1in}

Oriented rigid vertex graphs can have usual oriented nodes as well as special oriented nodes. 
By making different replacements for each kind of node we obtain more refined invariants for oriented RV graphs. For special oriented nodes we can choose to provide a marker as shown in Figure 18. We call this as special oriented vertex with indicator.

\begin{center}
     \begin{tabular}{c}
     \includegraphics[height=5cm]{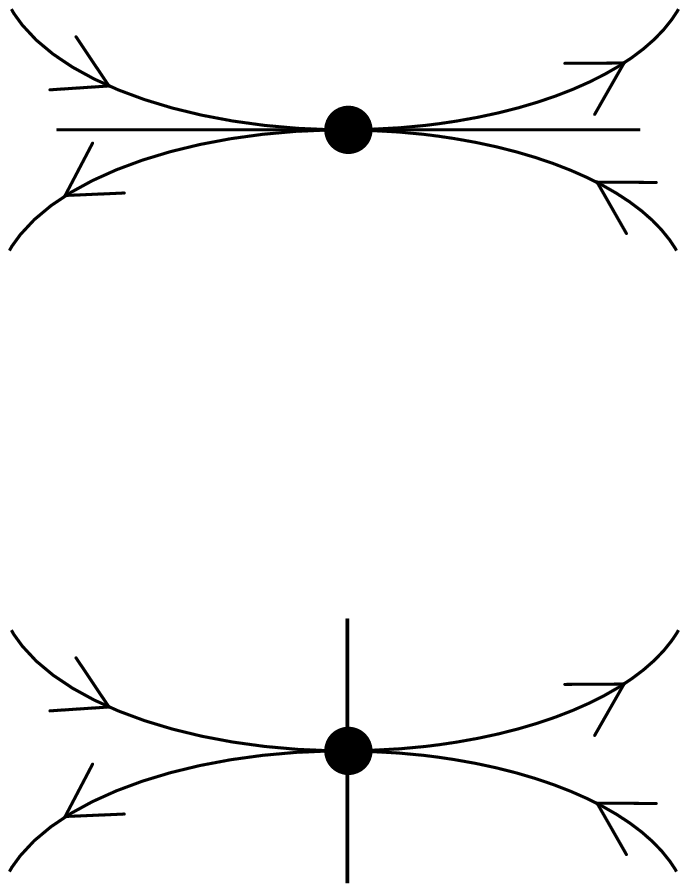}
     \end{tabular}
     \end{center}
     
    \centerline{\bf Figure 18 - Special Node Markers.}

Once the nodes have indicators then we can define the Gauss code by using the traverse at a crossing 
that corresponds to splitting it along the indicator. We will use this idea in the next section on RNA
folding.


     


\section{Associating a rigid vertex graph with special oriented node to RNA foldings} 
 
 It was discussed in \cite{Mag} that the topology of an RNA folding can be studied using its basic bonding node that arises from pairing of the base elements in the linear RNA chain.
This folding vertex  appears as shown in Fig 19. For the purpose of this discussion this folding vertex is
a special rigid vertex as shown in Figures 21 and 22. Thus we define an (abstract) protein folding
graph to be a graph with special vertices of this type. We often indicate the nodes as folding vertices
as shown in Figure 19. Showing the vertices in this way reminds us that that they can be ``unfolded" by 
cutting the node along the bonding lines. This is the same as cutting the node along the indicator for the special vertex in Figure 22.

\begin{center}
     \begin{tabular}{c}
     \includegraphics[height=3cm]{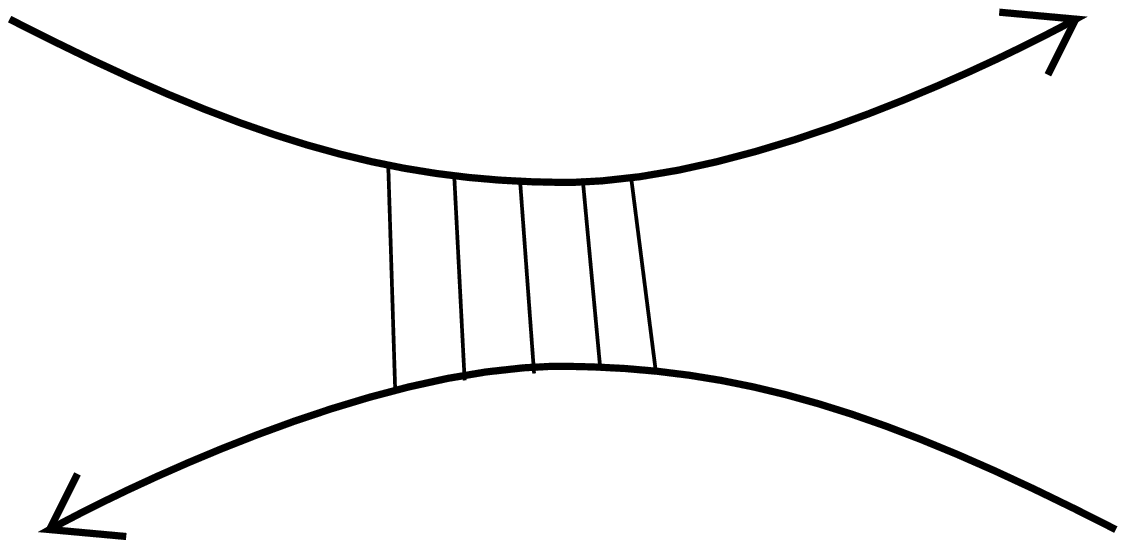}
     \end{tabular}
     \end{center}
     
    \centerline{\bf Figure 19 - Protein Folding Vertex.}
    
  \vspace{.1in}  

A folding vertex is assumed to be rigid, i.e., not subject to any twisting, while the oriented arcs that enter or leave the vertex are topologically flexible. Thus following moves shown in Fig 20 define isotopy for such a vertex.

 \begin{center}
     \begin{tabular}{c}
     \includegraphics[height=6cm]{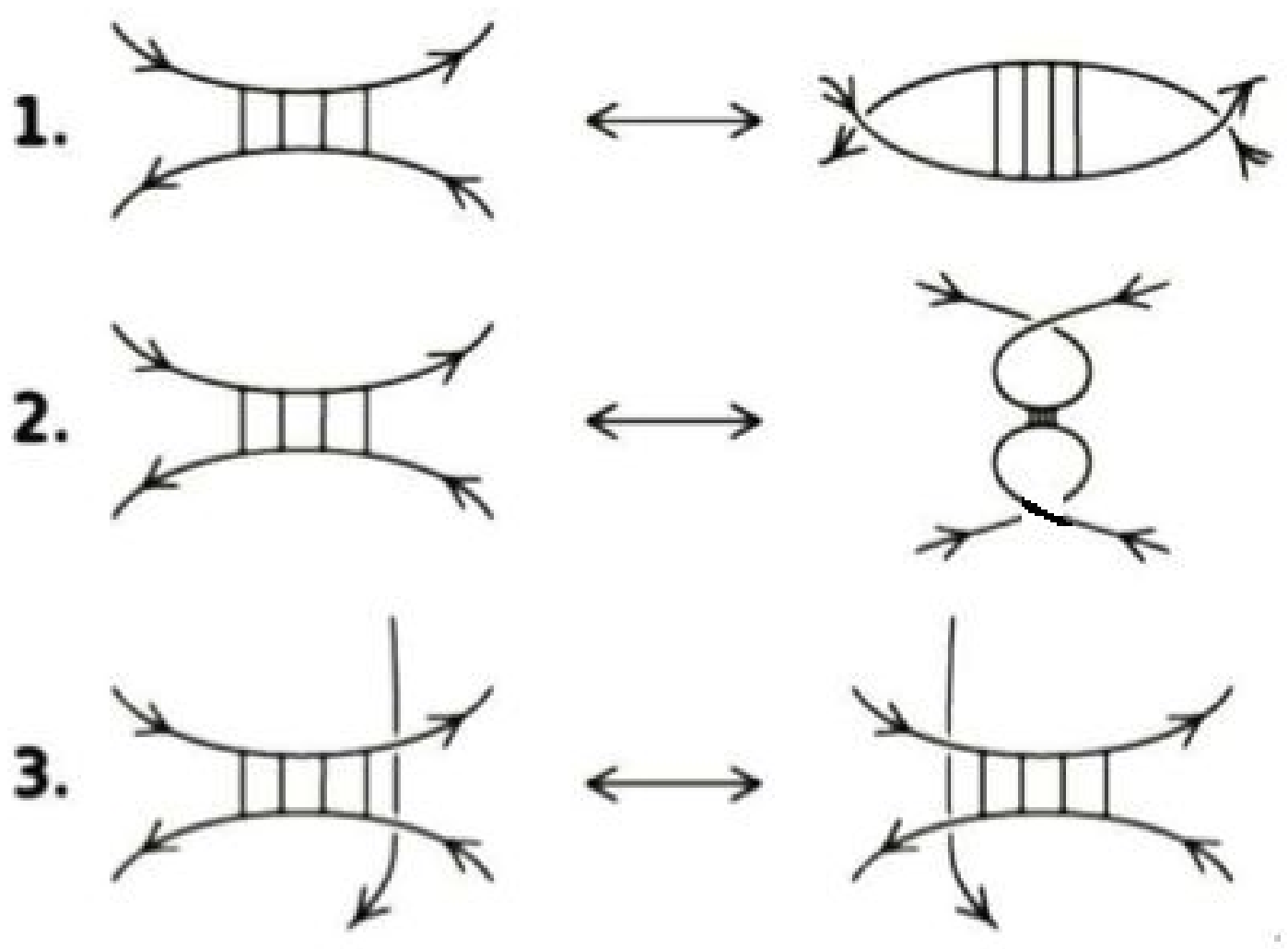}
     \end{tabular}
     \end{center}
     
    \centerline{\bf Figure 20 - Moves on Foldings.}
    
    \vspace{.1in}
 
 These moves along with the usual Reidemister Moves away from the bonds give us a complete characterization of rigid vertex isotopy.

We translate the folding vertex to a rigid special oriented vertex by the convention shown in Fig 21.
As a special oriented vertex it is identical to the translation in Figure 22.

\begin{center}
     \begin{tabular}{c}
     \includegraphics[height=2cm]{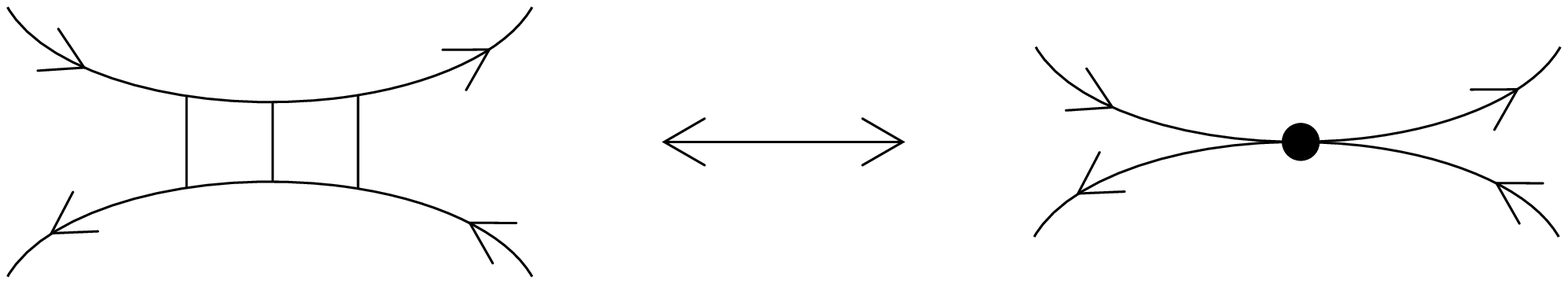}
     \end{tabular}
     \end{center}
     
    \centerline{\bf Figure 21 - Folding Vertex and Special Vertex.}

\begin{center}
     \begin{tabular}{c}
     \includegraphics[height=2cm]{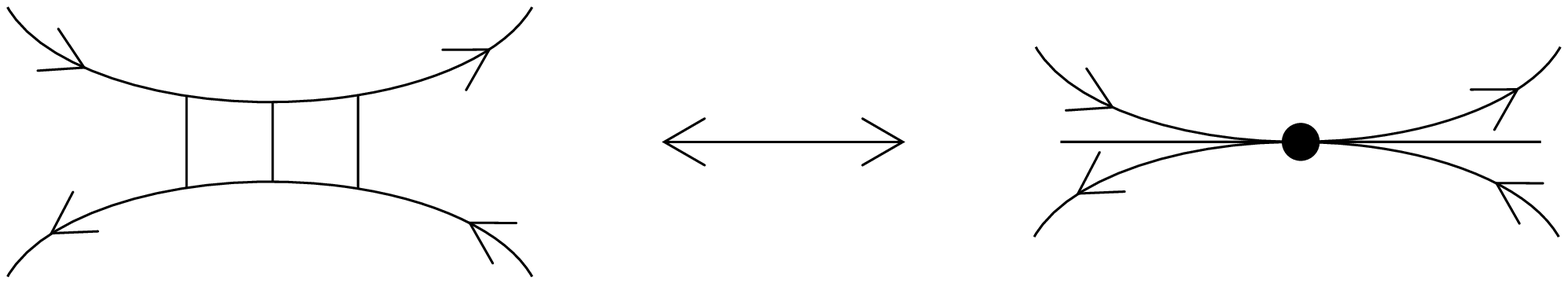}
     \end{tabular}
     \end{center}
     
    \centerline{\bf Figure 22 - Folding Vertex and Marked Special Vertex.}
    
\vspace{.1in}

Thus any invariant for a rigid vertex graph isotopy will serve as an invariant for RNA folding structure. We will see that a parity invariant turns out to be strong for classifying the RNA folding structure.

We discuss a few examples:

 Consider two foldings $T_1$ and $T_2$ as shown in Figure 23.

\begin{center}
     \begin{tabular}{c}
     \includegraphics[height=5cm]{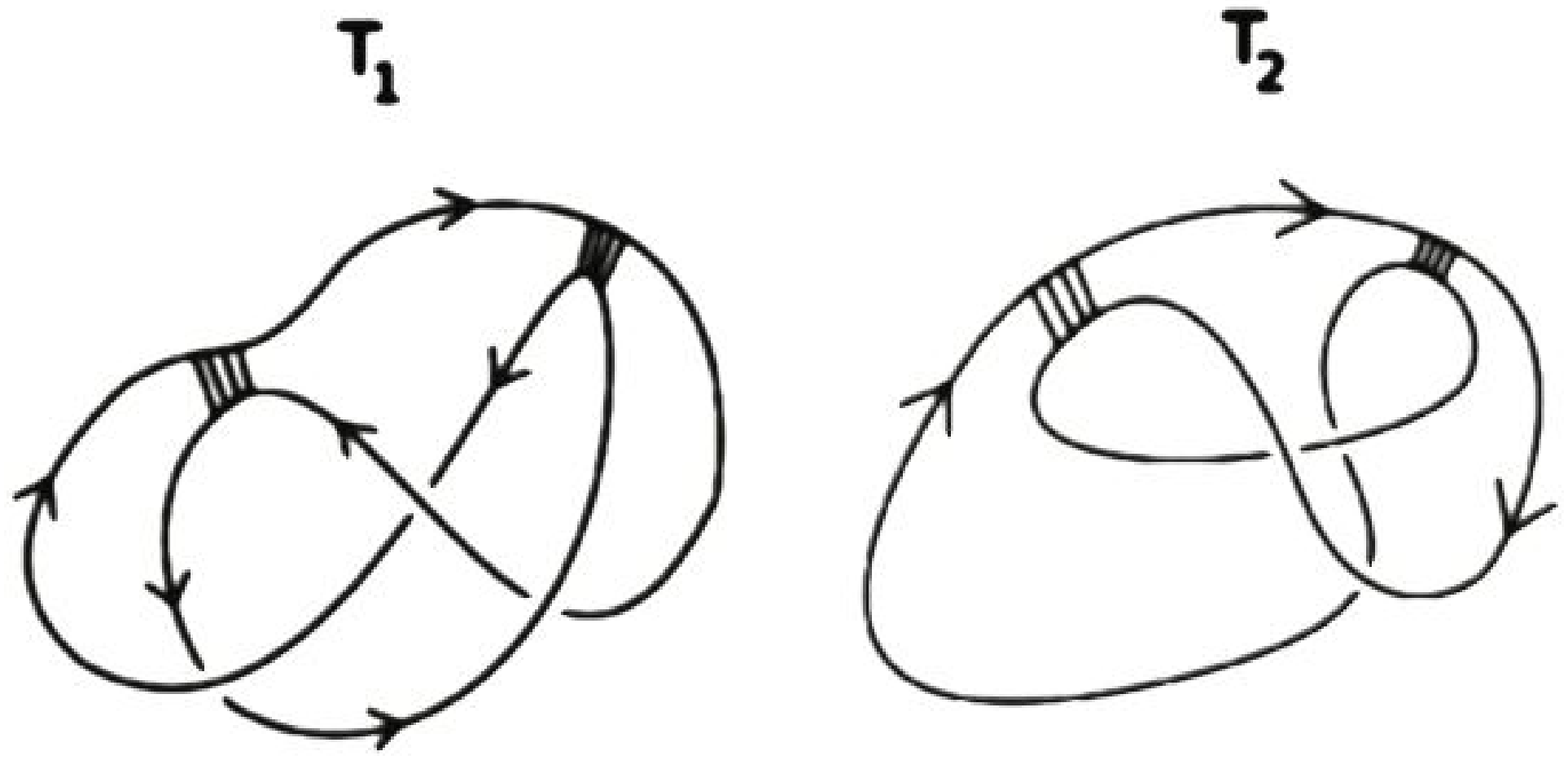}
     \end{tabular}
     \end{center}
     
    \centerline{\bf Figure 23 - Foldings.}
    
    \vspace{.1in}

We translate these foldings to rigid vertex graphs $G_1$ and $G_2$ respectively, with special oriented vertex, according to our convention discussed above. These graphs are shown in Figure 24.

\begin{center}
     \begin{tabular}{c}
     \includegraphics[height=4cm]{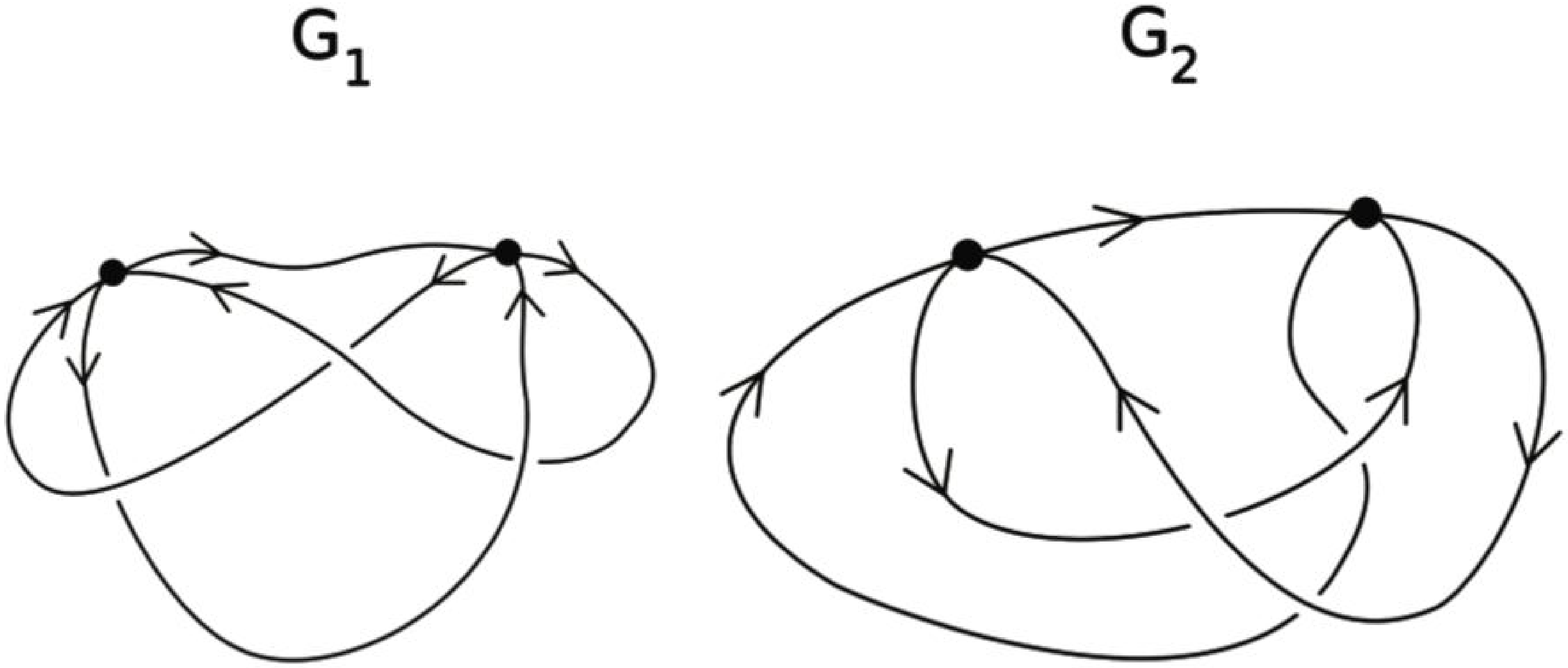}
     \end{tabular}
     \end{center}
     
    \centerline{\bf Figure 24 - Foldings Translated with Special Vertices.}
    
    \vspace{.1in}

In both these graphs, all nodes have odd parity. Note that in these graphs, with special vertices, there 
is only one smoothing at the vertex that preserves orientation. This smoothing corresponds to the
unfolding for the RNA folding vertex. Understanding this, one does not need a further indicator.
Hence we resolve them, preserving local orientation, and obtain knots $K_1$ and $K_2$ respectively as in Figure 25. In this case we see from the unfolding that $G_{1}$ is topologically different
from  $G_{2}.$
\smallbreak

In Figures 26 we give another example.
The corresponding rigid vertex graphs  can be unfolded, and it can be seen in Figure 26 that both of them result in unknots. Such foldings are known as {\em simple pseudoknots}. However these foldings are not equivalent. To ascertain that, we can replace the node by a full positive twist and then it is easy to see that the replacements are not isotopic.  We leave this as an exercise for the reader.

\begin{center}
     \begin{tabular}{c}
     \includegraphics[height=4cm]{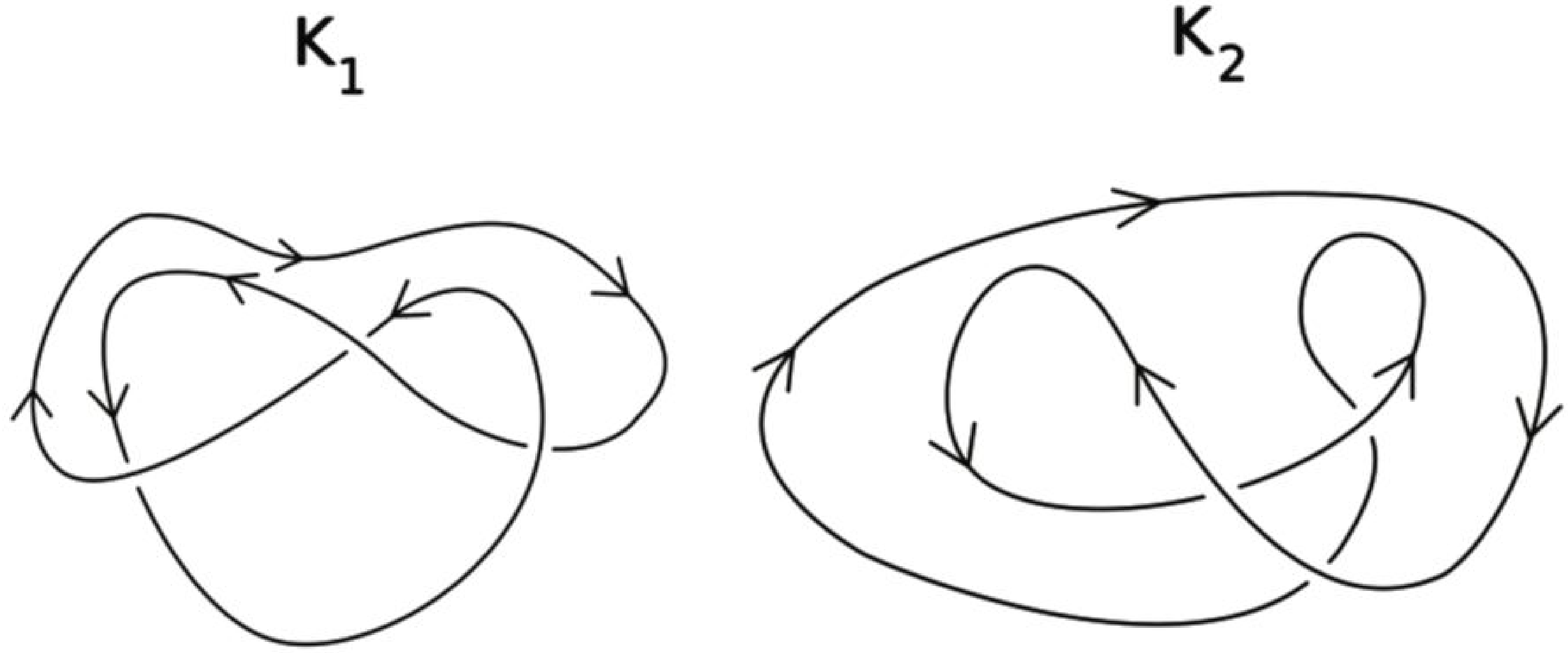}
     \end{tabular}
     \end{center}
     
    \centerline{\bf Figure 25 - Unfolding Special Vertices.}
    
    \vspace{.1in}

Let us take a more complex example of pseudo knots as shown in Figure 27. 
In $F_1$ all three nodes are of even parity, whereas in $F_2$ we have two nodes with odd parity and one node with even parity. Suppose we replace each even node by a full positive twist and an odd node by a full negative twist, we obtain knots $K_1$ and $K_2$ as shown in Figure 28. (Note that we are using the same names as in Figure 25. Let there be no confusion.)  It is easy to show that the knots $K_1$ and $K_2$ are not isotopic.

\begin{center}
     \begin{tabular}{c}
     \includegraphics[height=8cm]{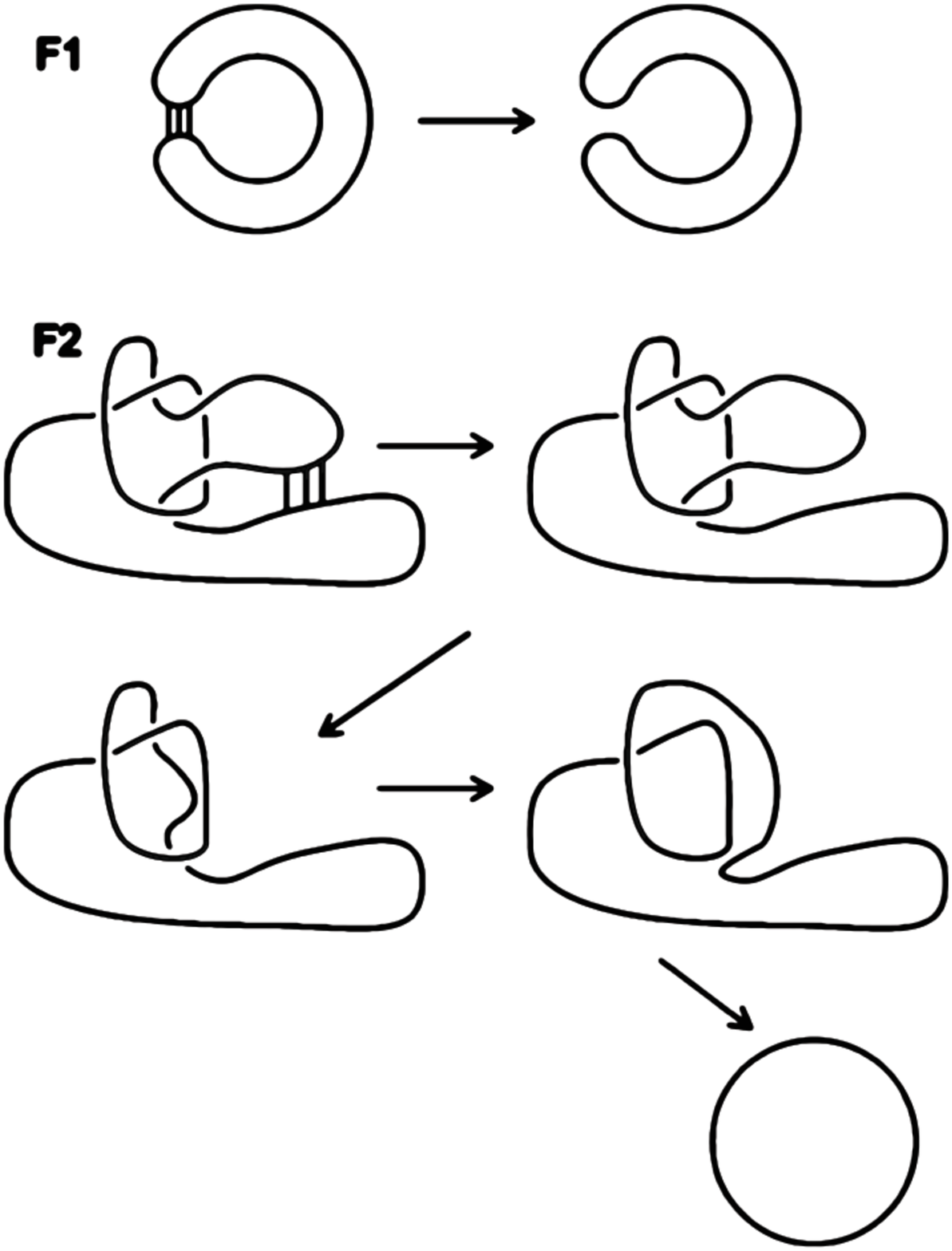}
     \end{tabular}
     \end{center}
     
    \centerline{\bf Figure 26 - Unfolding.} 
    
    \vspace{.1in}

     
    

\begin{center}
     \begin{tabular}{c}
     \includegraphics[height=5cm]{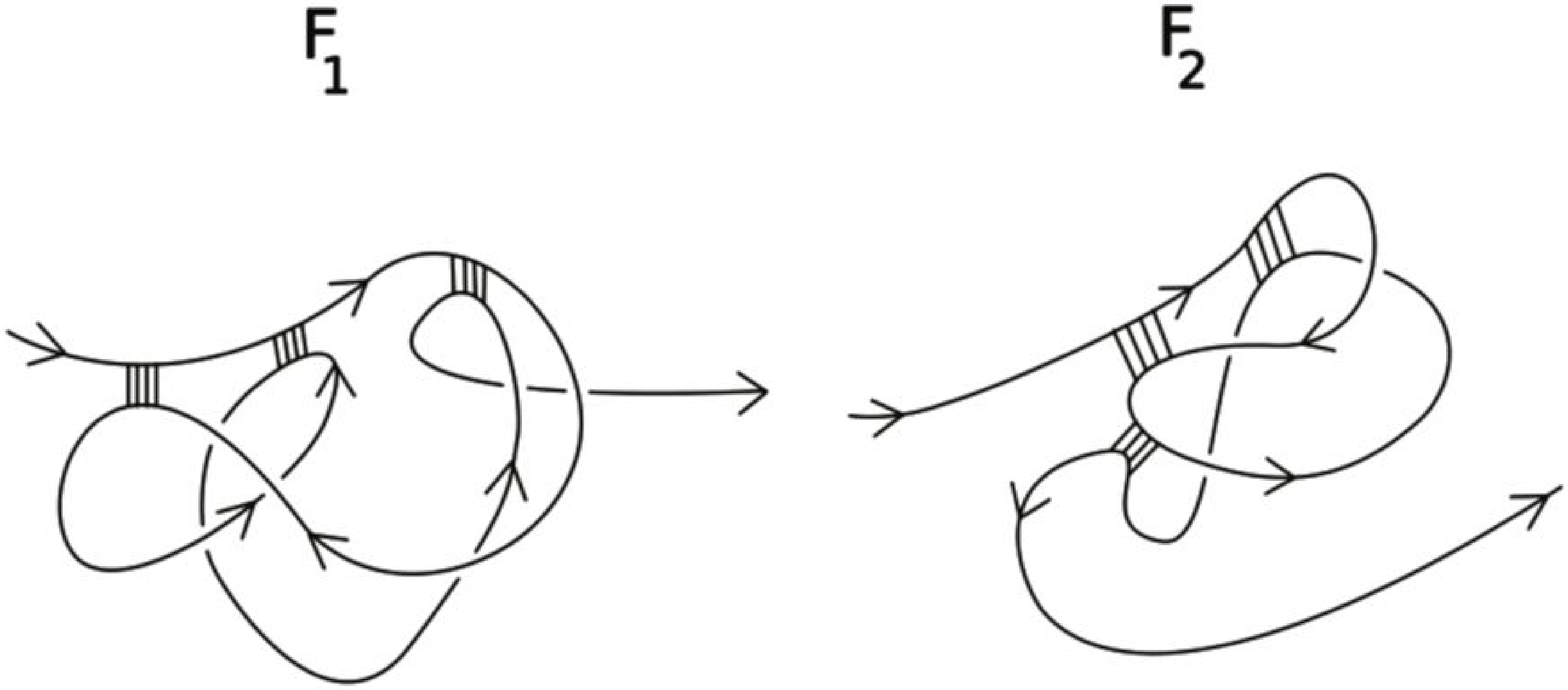}
     \end{tabular}
     \end{center}
     
    \centerline{\bf Figure 27 - Foldings.}
    
\vspace{.1in}

    \begin{center}
     \begin{tabular}{c}
     \includegraphics[height=5cm]{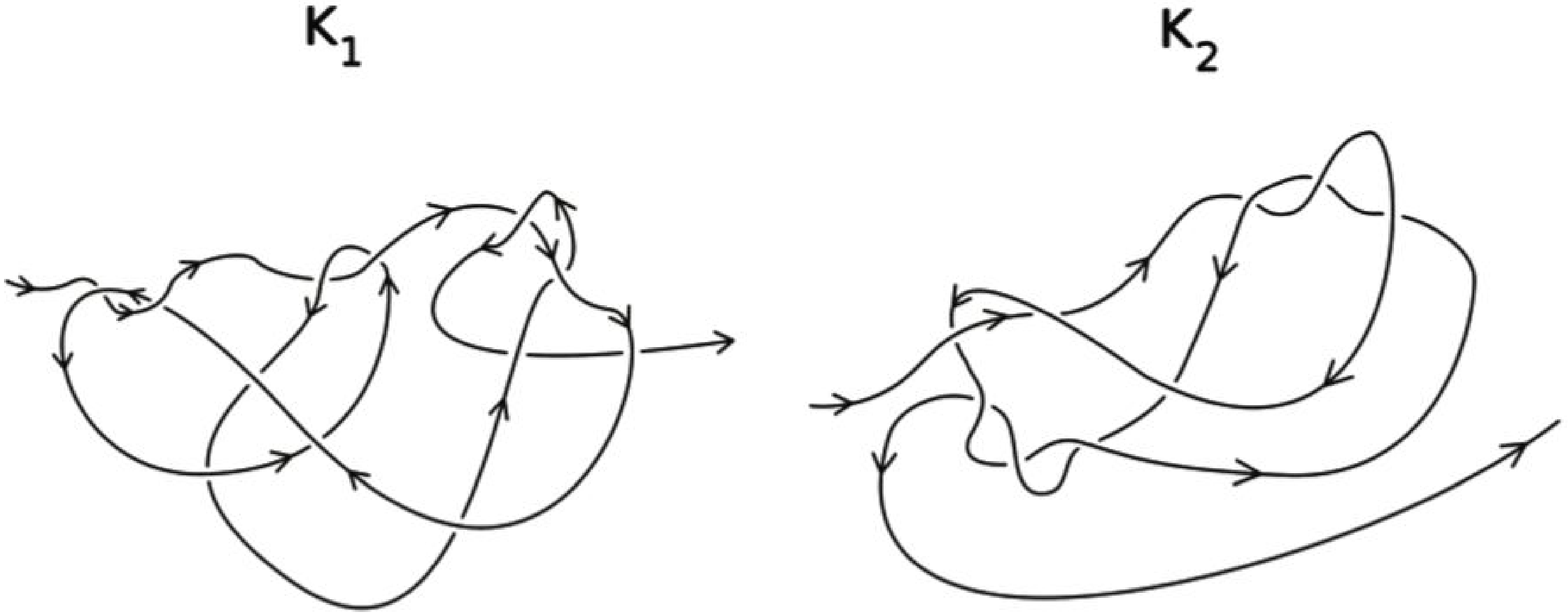}
     \end{tabular}
     \end{center}
     
    \centerline{\bf Figure 28 - Replacements.}
    
    \section{Virtual Knot Theory}
Knot theory
studies the embeddings of curves in three-dimensional space.  Virtual knot theory studies the  embeddings of curves in thickened surfaces of arbitrary
genus, up to the addition and removal of empty handles from the surface. Virtual knots have a special diagrammatic theory, described below,
that makes handling them
very similar to the handling of classical knot diagrams. Many structures in classical knot
theory generalize to the virtual domain.
\bigbreak  

In the diagrammatic theory of virtual knots one adds 
a {\em virtual crossing} (see Figure 29) that is neither an over-crossing
nor an under-crossing.  A virtual crossing is represented by two crossing segments with a small circle
placed around the crossing point. 
\bigbreak

Moves on virtual diagrams generalize the Reidemeister moves for classical knot and link
diagrams.  See Figure 29.  One can summarize the moves on virtual diagrams by saying that the classical crossings interact with
one another according to the usual Reidemeister moves while virtual crossings are artifacts of the attempt to draw the virtual structure in the plane. 
A segment of diagram consisting of a sequence of consecutive virtual crossings can be excised and a new connection made between the resulting
free ends. If the new connecting segment intersects the remaining diagram (transversally) then each new intersection is taken to be virtual.
Such an excision and reconnection is called a {\it detour move}.
Adding the global detour move to the Reidemeister moves completes the description of moves on virtual diagrams. In Figure 29 we illustrate a set of local
moves involving virtual crossings. The global detour move is
a consequence of  moves (B) and (C) in Figure 29. The detour move is illustrated in Figure 30.  Virtual knot and link diagrams that can be connected by a finite 
sequence of these moves are said to be {\it equivalent} or {\it virtually isotopic}.
\bigbreak

 \begin{center}
     \begin{tabular}{c}
     \includegraphics[height=7cm]{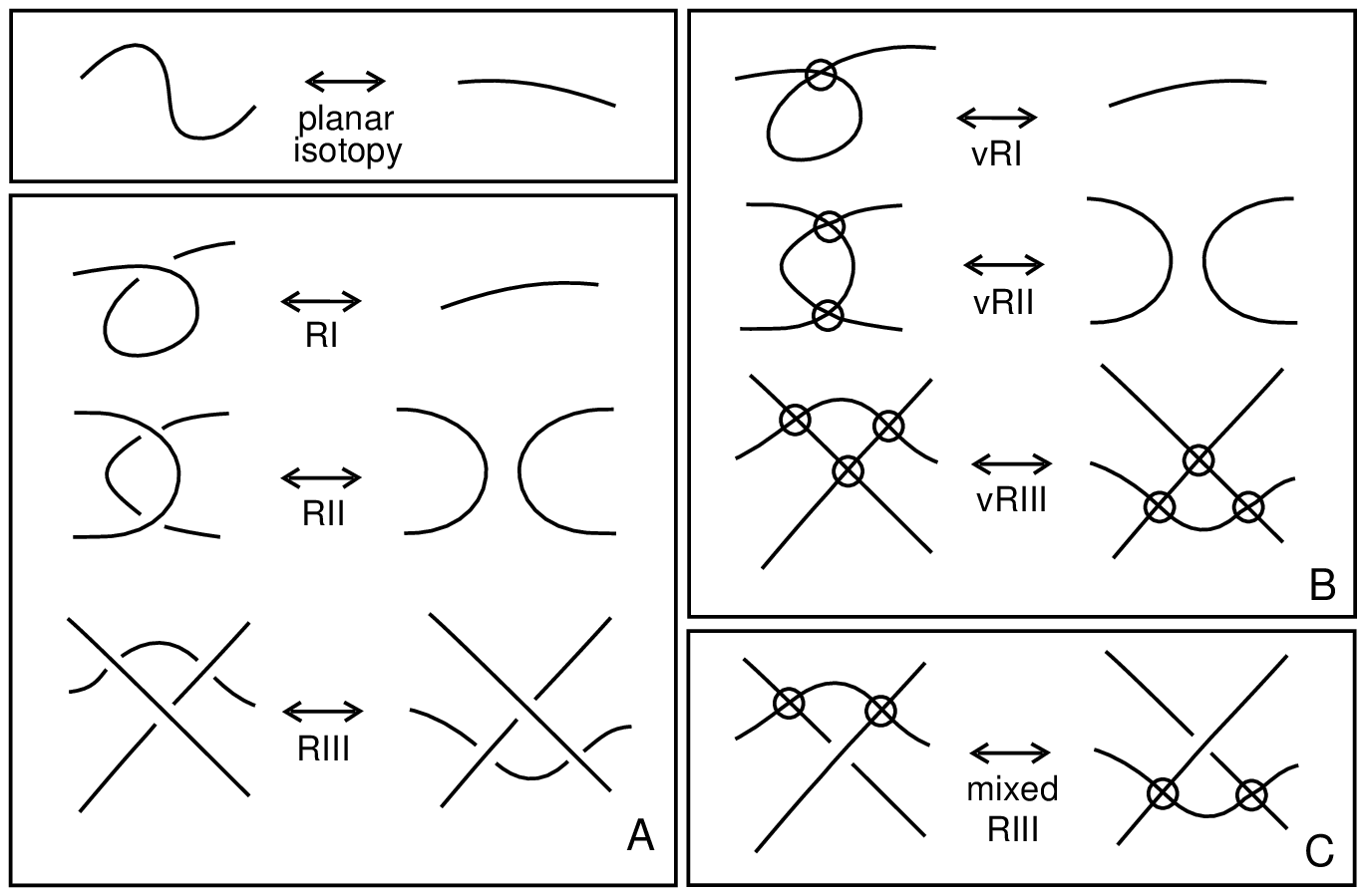}
     \end{tabular}
     \end{center}
     
    \centerline{\bf Figure 29 - Moves on Virtual Knots and Links.}

\begin{center}
     \begin{tabular}{c}
     \includegraphics[height=2cm]{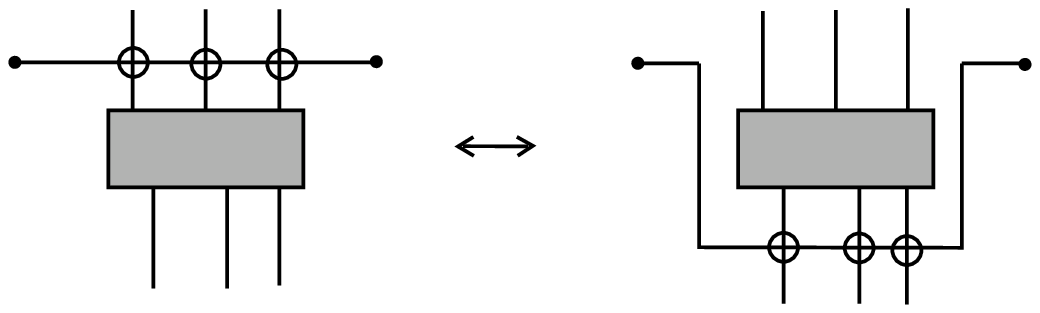}
     \end{tabular}
     \end{center}
     
    \centerline{\bf Figure 30 - Detour Move. }

\begin{center}
     \begin{tabular}{c}
     \includegraphics[height=2cm]{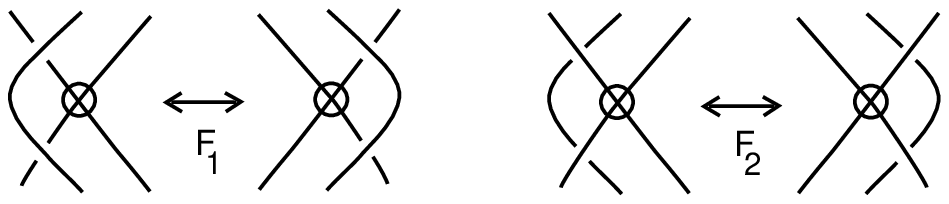}
     \end{tabular}
     \end{center}
     
    \centerline{\bf Figure 31 -  Forbidden Moves.}

\vspace{.1in}

Another way to understand virtual diagrams is to regard them as representatives for oriented Gauss codes \cite{GPV}, \cite{VKT,SVKT} 
(Gauss diagrams). Such codes do not always have planar realizations. An attempt to embed such a code in the plane
leads to the production of the virtual crossings. The detour move makes the particular choice of virtual crossings 
irrelevant. {\it Virtual isotopy is the same as the equivalence relation generated on the collection
of oriented Gauss codes by abstract Reidemeister moves on these codes.}  
\bigbreak

Figure 31 illustrates the two {\it forbidden moves}. Neither of these follows from Reidmeister moves plus detour move, and 
indeed it is not hard to construct examples of virtual knots that are non-trivial, but will become unknotted on the application of 
one or both of the forbidden moves. The forbidden moves change the structure of the Gauss code and, if desired, must be 
considered separately from the virtual knot theory proper. 
\bigbreak

\subsection{Interpretation of Virtuals Links as Stable Classes of Links in  Thickened Surfaces}
There is a useful topological interpretation \cite{VKT,DVK} for this virtual theory in terms of embeddings of links
in thickened surfaces.  Regard each 
virtual crossing as a shorthand for a detour of one of the arcs in the crossing through a 1-handle
that has been attached to the 2-sphere of the original diagram.  
By interpreting each virtual crossing in this way, we
obtain an embedding of a collection of circles into a thickened surface  $S_{g} \times R$ where $g$ is the 
number of virtual crossings in the original diagram $L$, $S_{g}$ is a compact oriented surface of genus $g$
and $R$ denotes the real line.  We say that two such surface embeddings are
{\em stably equivalent} if one can be obtained from another by isotopy in the thickened surfaces, 
homeomorphisms of the surfaces and the addition or subtraction of empty handles (i.e. the knot does not go through the handle).

\begin{center}
     \begin{tabular}{c}
     \includegraphics[height=4cm]{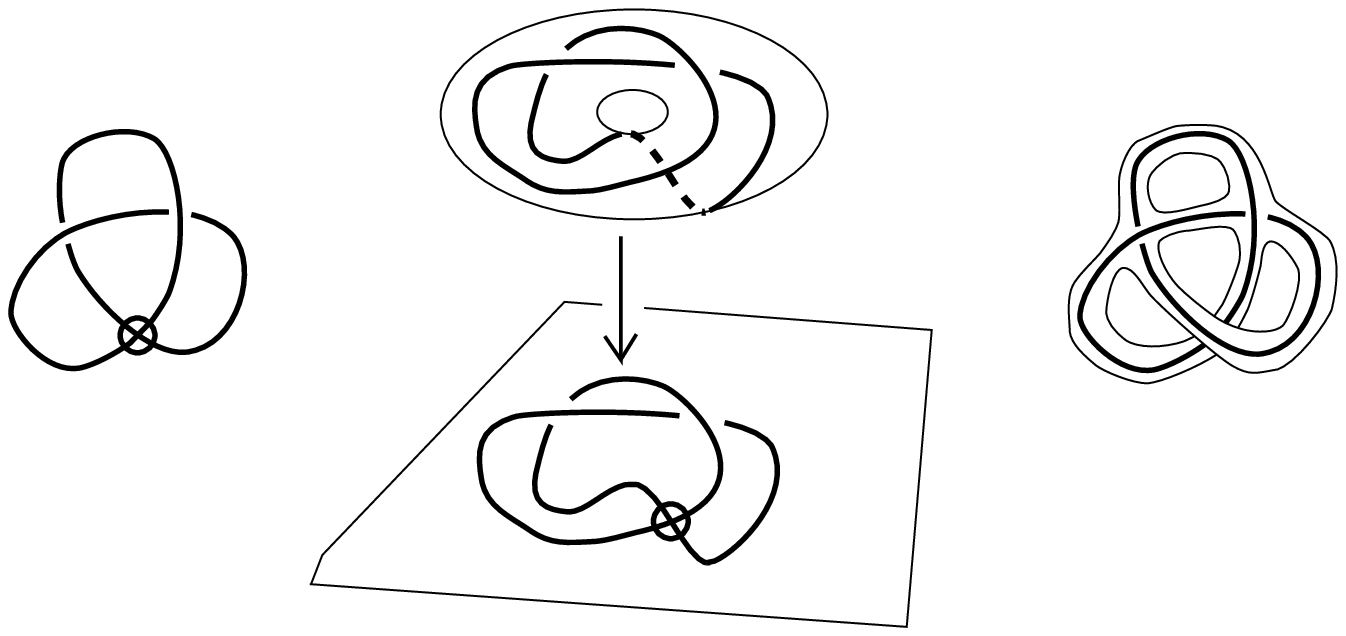}
     \end{tabular}
     \end{center}
     
    \centerline{\bf Figure 32 - Surfaces and Virtuals. }

\vspace{.1in}

\begin{center}
     \begin{tabular}{c}
     \includegraphics[height=4cm]{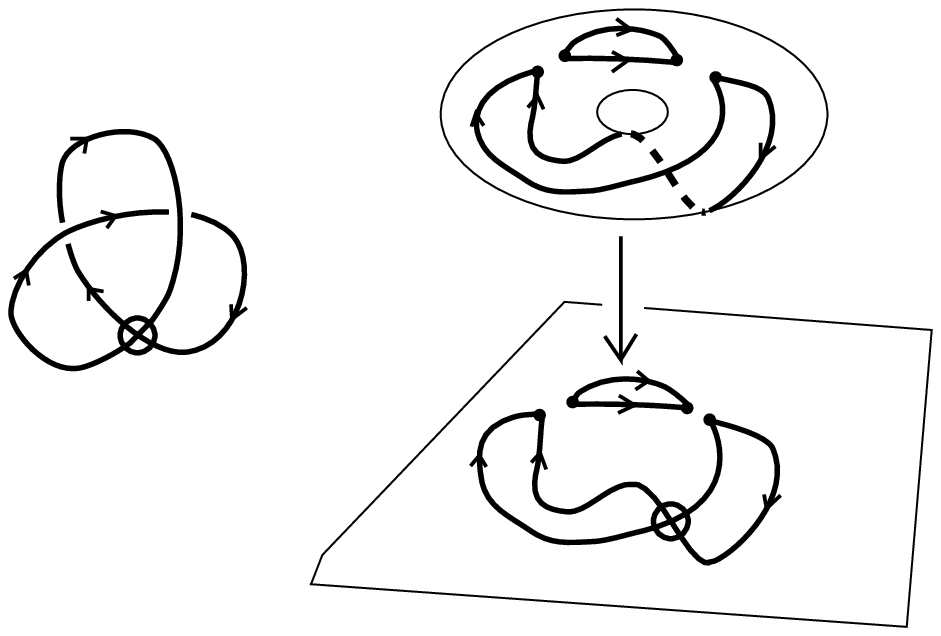}
     \end{tabular}
     \end{center}
     
    \centerline{\bf  Figure 33  - Surfaces and Virtual States. }

\noindent We have the
\smallbreak
\noindent
{\bf Theorem 1 \cite{VKT,DKT,DVK,Carter}.} {\em Two virtual link diagrams are isotopic if and only if their corresponding 
surface embeddings are stably equivalent.}  
\smallbreak
\noindent
\bigbreak  

\noindent In Figure 32 we illustrate some points about this association of virtual diagrams and knot and link diagrams on surfaces.
Note the projection of the knot diagram on the torus to a diagram in the plane (in the center of the figure) has a virtual crossing in the 
planar diagram where two arcs that do not form a crossing in the thickened surface project to the same point in the plane. In this way, virtual 
crossings can be regarded as artifacts of projection. The same figure shows a virtual diagram on the left and an ``abstract knot diagram" \cite{Kamada3,Carter} on the right.
The abstract knot diagram is a realization of the knot on the left in a thickened surface with boundary and it is obtained by making a neighborhood of the 
virtual diagram that resolves the virtual crossing into arcs that travel on separate bands. The virtual crossing appears as an artifact of the
projection of this surface to the plane. The reader will find more information about this correspondence \cite{VKT,DKT} in other papers by the author and in
the literature of virtual knot theory.
\bigbreak
 
 \subsection{Virtual graphs and virtual graph invariants}
 We now explain how to extend the results of our paper  to include graphs that have virtual crossings. We call such an entity a 
 {\em virtual graph} and by the same remarks as in the last section, we can interpret a virtual graph as a
graph embedding in a thickened surface, taken up to stabilization. In other words, the virtual knot theory extends directly to the study of such virtual graphs. One needs to extend the detour move so that virtual crossings can pass over graphical crossings. This is shown in Figure 34.  
\bigbreak

\begin{center}
     \begin{tabular}{c}
     \includegraphics[height=2cm]{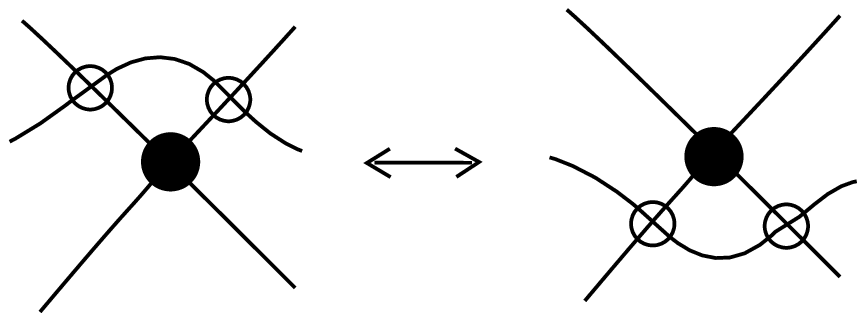}
     \end{tabular}
     \end{center}
     
    \centerline{\bf Figure 34 -  Detour Move for a Graphical Vertex. }

\bigbreak

We extend all of our definitions of parity, by simply ignoring the virtual crossings in making the corresponding parity counts. We are then prepared to apply all the methods of this paper to virtual
graph embeddings such as the example shown in Figure 35. In this example, the graph $G$ is seen to
be non-trivial by the following argument: The two graphical nodes have the same (odd) parity, so we
can resolve both of them positively as is shown in the diagram $R$ in Figure 35. We then see that
$R$ is a non-trivial virtual knot by applying a result \cite{VKT} in virtual knot theory that tells us that the 
virtual Jones polynomial of $R$ is the same as the virtual Jones polynomial of $S$ (in the same Figure 35) where $S$ is obtained from $R$ by smoothing the two virtual crossings that flank a classical crossing in $R.$ But $S$ is a trefoil knot and is known to have a non-trivial Jones polynomial.
This example illustrates how graph analysis as we have formulated it can contact virtual knot theory.
We leave the subject of virtual graphs at this point, but will return to it in a sequel to the present paper.
\bigbreak

\begin{center}
     \begin{tabular}{c}
     \includegraphics[height=6cm]{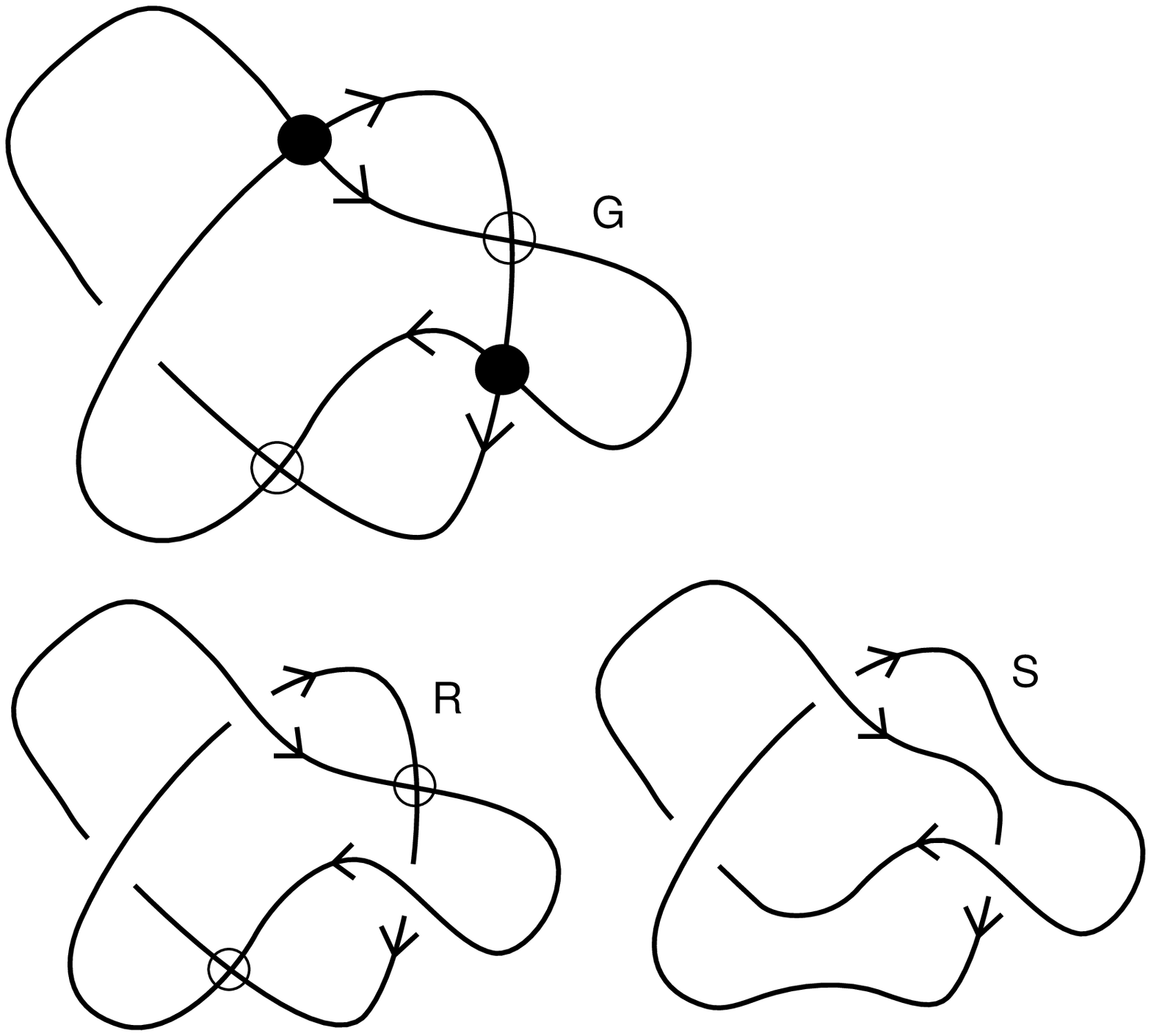}
     \end{tabular}
     \end{center}
     
    \centerline{\bf Figure 35  - Virtual Graph Example. }

\bigbreak

\noindent{\bf\large Acknowledgment:} Authors are thankful to two undergrad students Sarthak and Sriram for helping with the graphics.

\end{document}